%% file: Elements_Wrapper.tex
%
% Wrapper for Elements and introduction and afterword.
%

\def\wrapperUsed{WrapperPresent}
\magnification=\magstep1

\footline={\relax}

\
\vskip2.4truein

\font\TitleBBB=msbm10 scaled \magstep4
\centerline{\TitleBBB E~L~E~M~E~N~T~S}
\bigskip
\centerline{\bf Edward Nelson}

\bigskip
\bigskip
\centerline{With an}
\bigskip
\centerline{ Introduction by Sarah Jones Nelson}
\bigskip
\centerline{ and an }
\bigskip
\centerline{ Afterword by Sam Buss and Terence Tao }
\vfill
\vfill
\eject

\input SarahIntro_arXiv.tex

\vfill\eject

\input NelsonPhotoPage.tex

\pageno=1

\input elem.tex

\vfill\eject

\footline={\relax}\headline={\relax}

\input BussTao_Afterword_arXiv.tex

\end

%% file: SarahIntro_arXiv.tex
\font\sambigrm=cmr10 scaled \magstep2

{
\ifx\wrapperPresent\relax
\magnification=\magstep1
\fi
\hsize=5.0in
\hoffset=0.25in

\
\vskip1cm

\centerline{
\sambigrm  Introduction}

\vskip1cm

In September of 2011 Edward Nelson announced that he had a proof of the inconsistency of Peano Arithmetic. He had devoted twenty-five years to constructing the proof. When Terence Tao and Daniel Tausk independently found an error, Ed withdrew his claim at once and cheerfully returned to work on it the next day. By March of 2013 -- confident that he had corrected the error -- he wrote a project proposal with ``a crucial new insight: technically, bounds in the Hilbert-Ackermann consistency theorem depending only on rank and level, but not on the length of the proof.'' The project he proposed ``challenges the entire current understanding and practice of mathematics....It will radically change the way mathematics is done. This will affect the philosophy of mathematics, how the nature of mathematics is conceived. It will also affect the sciences that use mathematics, especially physics....It addresses a very big question: does mathematics consist in the discovery of truths about some uncreated eternal reality (the traditional Platonic view), or is it a humble human endeavor to construct abstract patterns that will be sound, free of all contradiction, in the hope that some will be beautiful, uplifting the human spirit, and that some (not necessarily different ones) will be of practical use to improve the human lot?''

The following excerpts from Ed's proposal describe his vision of a new mathematics and the open question of the consistency of Peano Arithmetic.
\vskip 1em

``Peano Arithmetic is one of the simplest and most fundamental of mathematical
theories. Its consistency, however, has not been proved by any means that all
mathematicians accept. It implies that all primitive recursive functions are
total, but they are directly defined only for numerals, and the argument that
the values always reduce to numerals is circular. The proposal is to complete a
proof that Peano Arithmetic is in fact inconsistent. The principal output will
be a book entitled `Elements'. The outcome will be a major change in the way
mathematics is done, with philosophical and scientific consequences.

``The guiding spirit of this investigation is that mathematics is not some
uncreated abstract reality that we can take for granted and explore, but that
it is a human endeavor in which one should begin by looking at the very
simplest concepts without taking them for granted.

``Numbers are constructed from $0$ by successively taking successors; ${\rm S
\ldots S}0$ is called a numeral. Definitions of primitive recursive functions,
such as addition, multiplication, exponentiation, superexponentiation, and so
forth, are schemata for constructing numerals. They define a value for $0$ and
then a value for ${\rm S}x$ in terms of the value for $x$. But when numerals
are substituted for the variables in such a schema it is not clear that it
defines a numeral: the putative number of steps needed to apply the definitions
can only be expressed in terms of the expressions themselves. The argument is a
vicious circle. Consequently, the consistency of Primitive Recursive
Arithmetic, and {\it a fortiori} of Peano Arithmetic ($\rm P$), is an open
question.

``Here is a nontechnical description of how I propose to show that $\rm P$ is
inconsistent. We start with a weaker theory $\rm Q$ and by relativization
techniques extend it to a stronger theory $\rm Q^*$. Proofs in $\rm Q^*$ reduce
to proofs in $\rm Q$. $\rm Q^*$ arithmetizes $\rm Q$ itself -- that is, it
expresses the syntax of $\rm Q$ by a term of $\rm Q^*$ that here I shall denote
by $\rm @Q$. Remarkably, $\rm Q^*$ proves that there is no open proof of a
contradiction in $\rm @Q$. (`Open' means that the proof has no quantifiers,
i.e., symbols for `there exists' and `for all'.) All of this was done in my
book `Predicative Arithmetic' (Princeton University Press, 1986) and is being
redone with complete proofs in the book `Elements', a work in progress that is
the subject of this proposal. The Hilbert-Ackermann consistency theorem implies
that there is no proof of a contradiction, even with quantifiers, in $\rm @Q$.
This theorem can be only partially established in $\rm Q^*$. The crucial new
insight is that it can be proved provided there are bounds on features of the
proof (called rank and level) but emphatically not depending on the length of
the proof. These bounds on rank and level cannot be proved in $\rm Q^*$, but for
each specific proof, $\rm P$ proves that $\rm Q^*$ proves them!  This opens the
way to exploit the stunning proof without self-reference of G\"{o}del's second
incompleteness theorem by Kritchman and Raz (Notices of the American
Mathematical Society, December 2011). The upshot is that $\rm P$ proves that
$\rm Q$ is inconsistent. But, as is well known, $\rm P$ also proves that $\rm
Q$ is consistent. Therefore Peano Arithmetic $\rm P$ is inconsistent.''
\vskip1em

Ed was the only living mathematician who could argue from purely syntactic reasoning without the traditional semantics established by Plato. John Conway suggested to me that this might explain why no one has fully understood Ed's deeply unique insights into the foundations of contemporary mathematics. He stood alone in the world, courageous as a formalist of a new ontology of integers: proof that completed infinities do not exist and that human minds invented numbers never discovered or revealed from platonic forms of any fundamental reality. I am hopeful that the mathematical community will boldly investigate ``Elements'' and the unshakeable foundations Ed sought to build.

\vskip1cm

\hfill Sarah Jones Nelson
\vskip 1ex

\hfill {\it September 10, 2015}

\vfill

}

\ifx\wrapperUsed\nil
\def\thisEnd{\end}
\else
\def\thisEnd{\relax}
\fi

\thisEnd

%% file: NelsonPhotoPage.tex
\ \vfill

\input epsf
\centerline{\epsfxsize=2in \epsfbox{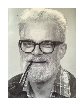}}
\vskip-6pt
\centerline{\hbox to 2in{\hfill \font\sixrm=cmr6 \sixrm Photo courtesy Sarah Jones Nelson}}
\vskip1em
\centerline{\font\twelverm=cmr12 \twelverm Edward Nelson}
\vfill

\eject

%% file: elem.tex
\magnification=\magstep1
\input NelsonarXivElemMacros.tex

\input hyperbasics.tex
\input colordvi.tex

\def\pdfklink#1#2{\hbox{\textBlue{#1}\textBlack}}

{\nopagenumbers \headline={} \footline={} \parindent=0pt

\phantom.
%\yy\pdfklink{http://world.std.com/~rjs/z138.pdf}{http://world.std.com/~rjs/z138.pdf}\zz

\vskip2.4truein

\centerline{\BBB E~L~E~M~E~N~T~S}

\bigskip
\bigskip
\bigskip
\bigskip

\centerline{\bf Edward Nelson}

\vfill\eject

\phantom.

\bigskip

{\ninerm Copyright \copyright\ 2013 by Edward Nelson, Department of Mathematics, Princeton University,
email: nelson@math.princeton.edu. All rights reserved.

\bigskip

This book may be reproduced in whole or in part for any noncommercial purpose.

\bigskip

All numbers in this book are fictitious, and any resemblance to actual numbers,
even or odd, is entirely coincidental.

\bigskip

This is a work in progress, version of March 12, 2013.

The latest version is posted at
\yy\pdfklink{http://www.math.princeton.edu/~nelson/books/elem.pdf}{http://www.math.princeton.edu/~nelson/books/elem.pdf}\zz

}

\vskip6.4truein

\centerline{$\varepsilon_0\ \omega^2\ 10\ 9\ 8\ 7\ 6\ 5\ 4\ 3\ 2\ 1\ 0$}

\vfill\eject

\phantom.

\vskip3.5truein

\centerline{To Sarah}

\medskip

\centerline{\it because you are overflowing with fun and faith}

\vfill\eject

\parindent=.6truein

\phantom.

\vskip2.5truein

Glory be to God for dappled things---

\quad For skies of couple-color as a brinded cow;

\qquad For rose-moles all in stipple upon trout that swim;

Fresh-firecoal chestnut-falls; finches' wings;

\quad Landscape plotted and pieced---fold, fallow, and plough;

\qquad And \acute all trades, their gear and tackle and trim.

\medskip

All things counter, original, sp\acute are, strange;

\quad Whatever is fickle, freckl\`ed (who knows how?)

\qquad With sw\acute\i ft, sl\acute ow; sweet, s\acute our; ad\acute azzle, d\acute\i m;

He fathers-forth whose beauty is p\acute ast change:

\medskip

\qquad\qquad\qquad\qquad\qquad\qquad\qquad\quad Pr\acute aise h\acute\i m.

\bigskip

\hfill\hbox{--- \ninerm G{\sevenrm ERARD} M{\sevenrm ANLEY}~H{\sevenrm OPKINS}}
\phantom{xxxxxxxxxx}

\vfill\eject

\phantom.

\vskip2truein

\centerline{\bigrm PART I}

\bigskip
\bigskip
\bigskip

\centerline{INCONSISTENCY}

\vfill\eject

\phantom.

\vfill\eject

} % end nopagenumbers etc.

\chap1.{Thesis}

The aim of this work is to show that contemporary mathematics, including
Peano Arithmetic, is inconsistent, to construct firm foundations for mathematics,
and to begin building on those foundations.

\section//   1.{Potential versus actual infinity}
Let us distinguish between the \ii concrete.\/ (i.e., the genetic in the
sense of pertaining to
origins, the given, what is present) on the one hand, and on the other the \ii formal.\/
(i.e., the abstract, the hypothetical).
All of mathematical activity is concrete, though the subject matter is formal.

A \ii numeral.\/ is a variable-free term of the language whose nonlogical symbols are
the constant~0 (\ii zero.\/) and the unary function symbol~S (\ii successor.\/).
Thus the numerals are 0,~S0, SS0, SSS0, \dots.

Numerals constitute a \ii potential infinity.. Given any numeral, we can construct a new
numeral by prefixing it with~S. Now imagine this potential infinity to be completed.
Imagine the inexhaustible process of constructing numerals somehow to have been finished,
and call the result \ii the set of all numbers., denoted by $\omega$ or \N.
Thus $\omega$ is thought to be an \ii actual infinity.\/ or a
\ii completed infinity.. This is curious terminology,
since the etymology of ``infinite'' is ``not finished''.

As a concrete concept, the notion of numeral is clear.
The attempt to formalize the concept usually proceeds as follows:

\nok// 1. zero is a number

\nok// 2. the successor of a number is a number

\nok// 3. zero is not the successor of any number

\nok// 4. different numbers have different successors

\nok// 5. something is a number only if it is so by virtue of $(/ 1)$ and $(/ 2)$

\medskip

We shall refer to this as the \ii usual definition..
Sometimes $(/ 3)$ and~$(/ 4)$ are not stated explicitly, but
it is the extremal clause~$(/ 5)$ that is unclear.
What is the meaning of ``by virtue of''?
It is obviously circular to define a number as something constructible by
applying $(/ 1)$ and~$(/ 2)$ any number of times. We cannot characterize numbers from
below, so we attempt to characterize them from above.

The study of the foundations of arithmetic began in earnest with
Dedekind\foot{Richard Dedekind, {\it Was sind und was sollen die Zahlen?},
Vieweg, Braunschweig, (1888).
English translation as {\it The nature and meaning of numbers} by W. W. Beman in
Richard Dedekind,
Essays on the Theory of Numbers, Open Court Publishing Company,
LaSalle, Illinois. Reprinted
by Dover, New York, 1963.} and
Peano.\foot{Giuseppe Peano, {\it Arithmetices principia, nova methodo exposita},
Bocca, Turin, (1889).
English translation in From Frege to G\umlaut odel: A source book in
mathematical logic, 1879--1931, ed. J. van Heijenoort,
Harvard University Press, Cambridge, Massachusetts, 1967.}
Both of these authors gave what today would be called set-theoretic
foundations for arithmetic.
In \ZFC\ (Zermelo-Fraenkel set theory with the axiom of choice), let us write
0 for the empty set and define the successor by

\noj $\'S'x=x\cup\lbrace x\rbrace $

\noin where the unary function symbol $\{\ \}$ (\ii singleton.\/)
and the binary function symbol~$\cup$ (\ii union.\/) have been defined
in the usual way. We define

\noj $X \' is inductive' \iff 0\in X \and \forall x[x\in X \impP \'S'x\in X]$

\noin Then the axiom of infinity of \ZFC\ is

\noj $\exists X[X \' is inductive']$

\noin and one easily proves in \ZFC\ that there exists a
unique smallest inductive set; i.e.,

\noj $\exists!X\{X \' is inductive' \and \forall Y[Y \' is inductive'
\impP X\subseteq Y]\}$

\noin where the binary predicate symbol $\subseteq$ has been defined as usual.
We define the constant~$\omega$ to be this smallest
inductive set:

\noj $\omega=X \iff X \' is inductive' \and
\forall Y[Y \' is inductive' \impP X\subseteq Y]$

\noin and we define

\noj $x \' is a number' \iff x \in \omega$
\medskip

Then the following are theorems:

\nok// 6. $0 \' is a number'$

\nok// 7. $x \' is a number' \imp \'S'x \' is a number'$

\nok// 8. $x \' is a number' \imp \'S'x\ne0$

\nok// 9. $x \' is a number' \and y \' is a number' \and x\ne y \imp
\'S'x\ne\'S'y$
\medskip

These theorems are a direct expression of $(/ 1)$--$(/ 4)$ of the usual definition.
But can we express the extremal clause~$(/ 5)$? The induction theorem

\noj $x \' is a number' \and Y \' is inductive' \imp x \in Y$

\noin merely asserts that for any property that can be expressed by a set in \ZFC, if
0~has the property, and if the successor of every element that has the property also
has the property, then every number has the property.

We cannot say, ``For all numbers~$x$
there exists a numeral~d such that $x=\ro d$'' since this is a category mistake
conflating the formal with the concrete.

Using all the power of modern mathematics, let us try to formalize the concept
of number.
Let T be any theory whose language contains the constant 0, the unary function
symbol~S, and the unary predicate symbol ``is a number'', such that $(/ 6)$--$(/ 9)$
are theorems of~T. For example, T~could be the extension by definitions of \ZFC\
described above or it could be \PA\ (Peano Arithmetic; see below) with the definition:
$x \' is a number' \iffF x=x$.

Have we captured the intended meaning of the extremal clause $(/ 5)$? To study this
question,
construct $\ro T^\phi$ by adjoining a new unary predicate symbol $\phi$ and the
axioms

\nok// 10. $\phi(0)$

\nok// 11. $\phi(x) \impP \phi(\'S'x)$

\noin Notice that $\phi$ is an undefined symbol. If T is \ZFC, we cannot form
the set $\lbrace x\in\omega\,\colon\phi(x)\rbrace$
because the subset axioms of \ZFC\ refer only to
formulas of \ZFC\ and $\phi(x)$ is not such a formula. Sets are not concrete objects,
and to ask whether a set with a certain property exists is to ask whether a certain
formula beginning with~$\exists$ can be proved in the theory.
Similarly, if T is \PA\ we cannot apply induction to $\phi(x)$ since this
is not a formula of \PA. Induction is not a truth; it is an axiom schema of a formal
theory.

If T is consistent then so is $\ro T^\phi$, because we can interpret
$\phi(x)$ by $x=x$. (And conversely, of course, if T~is inconsistent then
so is~$\ro T^\phi$.)
For any concretely given
numeral S\dots S0 we can prove $\rm\phi(S\ldots S0)$ in as many steps as there are
occurrences of~S in S\dots S0,
using these two axioms and modus ponens.

Proving~$\rm\phi(d)$ in~$\ro T^\phi$ perfectly expresses the notion that d~is
a number by virtue of $(/ 1)$ and~$(/ 2)$ of the usual definition.
We can read $\phi(x)$ as ``$x$~is a number by virtue of $(/ 1)$ and~$(/ 2)$''.
Therefore we ask: can

\nok// 12. $x \' is a number' \imp \phi(x)$

\noin be proved in $\ro T^\phi$? (That is, can we prove that our formalization
``$x$~is a number'' captures the intended meaning of the extremal clause?)
Trivially yes if T is inconsistent, so assume
that T~is consistent. Then the answer is no. Here is a semantic argument for
this assertion.

By $(/ 6)$--$(/ 9)$, none of the formulas

\noj $x \' is a number' \imp x=0 \,\orR\, x=\'S'0 \,\orR\,\cdots\,\orR\,
x=\'S'\ldots\'S'0$

\noin is a theorem of T. Hence the theory $\ro T_1$ obtained
from~T by adjoining a new constant~e and the axioms

\noj $\ro e \' is a number'$, $\ro e\ne0$, $\ro e\ne\'S'0$, \dots,
$\rm e\ne S\ldots S0$, \dots

\noin is consistent. By the G\umlaut odel completeness
theorem [G\umlaut o29]\foot{Kurt G\umlaut odel, {\it \umlaut Uber die
Vollst\umlaut andigkeit des
Logikkalk\umlaut uls}, doctoral dissertation, University of Vienna, (1929).}
[Sh~\S4.2],\foot{Joseph R. Shoenfield, {\it Mathematical Logic}, Association for
Symbolic Logic, A K Peters, Ltd., Natick, Massachusetts, 1967.}
$\ro T_1$~has a model $\sigma$ [Sh~\S2.6]. Let I be the
smallest subset of the universe~U
of the model containing $\sigma(0)$ and closed under
the function $\rm\sigma(S)$. Then $\rm\sigma(e)$ is not
in~I. Expand $\sigma$ to be a model $\sigma^\phi$ of~$\ro T^\phi$ by
letting~$\sigma^\phi(\phi)$ be~I. Then $(/ 12)$ is not valid in this model,
and so is not a theorem of~$\ro T^\phi$.

The conclusion to be drawn from this argument is that it is impossible to formalize
the notion of number in such a way that the extremal clause holds.%
\yy%
\foot{\pdfklink{http://frankandernest.com/cgi/view/display.pl?104-12-15}%
{http://frankandernest.com/cgi/view/display.pl?104-12-15}}
\zz

Most uses of infinity in mainstream mathematics involve infinity as a limit of
finite objects. The situation is qualitatively different when it comes to
the notion of truth in arithmetic. Those who believe that Peano Arithmetic {\it must\/}
be consistent because it has a model with universe~$\omega$ are placing their faith in a
notion of truth in that model that is utterly beyond any finite computation.

Despite all the accumulated evidence to the contrary,
mathematicians persist in believing in~$\omega$ as a real object
existing independently of any formal human construction. In a way this is not
surprising. Mathematics as a deductive discipline was
invented by Pythagoras, possibly with some influence from Thales.
The Pythagorean religion held that all is number, that the numbers
are pre-existing and independent of human thought. Plato was strongly influenced
by Pythagoras and has been called the greatest of the Pythagoreans.
Over two and a half millennia after Pythagoras,
most mathematicians continue to hold a religious belief in~$\omega$
as an object existing independently of formal human construction.

\bigskip

\hfill\vbox{\baselineskip=10pt\halign{\ninerm#\cr
Infinite totalities do not exist in any sense of the word\hfil\cr
(i.e., either really or ideally). More precisely, any mention,\hfil\cr
or purported mention, of infinite totalities is, literally,\hfil\cr
{\nineit meaningless.}\hfil\cr
}}\smallskip

\hfill\hbox{--- \ninerm A{\sevenrm BRAHAM} \ninerm R{\sevenrm OBINSON}}

\section//   2.{Against finitism}
There is obviously something inelegant about making arithmetic depend
on set theory. What today is called \ii Peano Arithmetic.\/ (\PA) is the
theory whose nonlogical symbols are the constant~0, the unary function symbol~S,
and the binary function symbols~$+$ and~$\.$, and
whose nonlogical axioms are

\nok// 13. $\'S'x\ne0$

\nok// 14. $\'S'x=\'S'y\imp x=y$

\nok// 15. $x+0=x$

\nok// 16. $x+\'S'y=\'S'(x+y)$

\nok// 17. $x\.0=0$

\nok// 18. $x\.\'S'y=(x\.y)+x$

\noin and all \ii induction formulas.

\nok// 19. $\rm A_x(0) \and \forall x[A \impP A_x(\'S'\ro x)] \imp \ro A$

\noin where A is any formula in the language of \PA.

The induction axiom schema $(/ 19)$ is usually justified as follows.
Assume the \ii basis.\/ $\rm A_x(0)$
and the \ii induction step.\/ $\rm \forall x[A\impP A_x(\'S'\ro x)]$.
For any numeral~d, a special case of the induction step is $\rm A_x(d)\impP A_x(Sd)$.
Then we can prove $\rm A_x(SSS0)$, starting with the basis and using
modus ponens three times, and similarly we can prove $\rm A_x(d)$ for any
numeral~d in as many steps as there are occurrences of~S in~d.

This justification leaves much to be desired. The argument applies only to numerals,
and shows that there is no need at all to postulate induction for numerals. By an
unspoken conflation of the genetic concept of numeral with the formal concept of
number, induction is postulated for numbers.

Induction is justified by appeal to the \ii finitary credo.\/: for every number~$x$ there
exists a numeral~d such that $x$ is~d.
It is necessary to make this precise; as someone once said,
it depends on what you mean by ``is''.
We cannot express it as a formula of arithmetic because ``there exists'' in
``there exists a numeral~d'' is a metamathematical existence assertion, not an
arithmetical formula beginning with~$\exists$.

Let A be a formula of~\PA\ with just one free variable~$x$ such that
$\vdash_{\sPA}\exists x\ro A$. By the least number principle (which is a form of
induction), there is a least number~$x$ satisfying~A. One might, thinking that every
number is a numeral, assert that there exists a numeral~d such that
$\vdash_{\sPA}\ro A_x(\ro d)$. But, as I learned from Simon Kochen, this does not work.
Let B be the formula asserting that there is no arithmetized proof of a contradiction
all of whose formulas are of rank at most~$x$. Then each $\ro B_x(\ro d)$ can be
proved in~\PA\ by introducing a truth definition for formulas of rank at most~d.
But if \PA\ is consistent then $\forall x\ro B$ is not a theorem of \PA, by
G\umlaut odel's second theorem. Now let A~be $\ro B\impP\forall x\ro B$. Then
$\exists x\ro A$ is a theorem of~\PA\ since it is equivalent to the tautology
$\forall x\ro B\impP\forall x\ro B$, but (if \PA~is consistent) there is no numeral~d
such that $\vdash_{\sPA}\ro A_x(\ro d)$.

The finitary credo can be formulated precisely using the concept of the standard model of
arithmetic: for every element~$\xi$ of~\N\ there exists a numeral~d such that it can be
proved that d~is equal to the name of~$\xi$; but this brings us back to set theory.
The finitary credo has an infinitary foundation.

The use of induction goes far beyond the application to numerals.
It is used to create new kinds of numbers (exponential, superexponential, and
so forth) in the belief that they already exist in a completed infinity.
If there were a completed infinity $\N$ consisting of all numbers, then the
axioms of \PA\ would be true assertions about numbers and \PA\ would be consistent.

It is not a priori obvious that \PA\ can express combinatorics, but
this is well known thanks to G\umlaut odel's great
paper [G\umlaut o31]\foot{Kurt G\umlaut odel, {\it \umlaut Uber formal unentscheidbare
S\umlaut atze der
Principia Mathematica und verwanter Systeme}, Monatshefte f\umlaut ur Math.
und Physik, 38 (1931) 173-198. English translation in From Frege
to G\umlaut odel: A source book in mathematical logic, 1879--1931, 596-616,
ed. Jean van Heijenoort, Harvard University Press, 1967.}
on incompleteness. As a consequence, exponentiation~$\uparrow$
and superexponentiation~$\Uparrow$ can be defined in \PA\ so that we have

\nok// 20. $x\uparrow 0=\'S'0$

\nok// 21. $x\uparrow\'S'y=x\.(x\uparrow y)$

\nok// 22. $x\Uparrow0=\'S'0$

\nok// 23. $x\Uparrow\'S'y=x\uparrow(x\Uparrow y)$

\noin and similarly for primitive-recursive functions in general.

Now substitute variable-free terms a and b for $x$ and $y$ in these equations. The left hand sides
construct, for example, $\rm a\Uparrow b$ when b is a numeral: first $\rm a\Uparrow0$, then
$\rm a\Uparrow S0$, $\rm a\Uparrow SS0$, and so on. But there is a problem: we write $\rm a\Uparrow Sb=
a\uparrow(a\Uparrow b)$, but $\rm a\uparrow c$ is defined only for numerals~c, and $\rm a\Uparrow b$
is not a numeral.

The finitary credo asserts that applications of these equations a sufficient number of times reduce
variable-free terms to numerals. But what evidence is there for this belief? The number of steps required
can only be expressed in terms of these primitive-recursive terms themselves---the reasoning is circular.
The objection being raised here is not some vague ``ultrafinitistic'' semantic argument that these numbers
are so large they don't really exist; the problem is structural. Indeed, to call such variable-free terms
``numbers'' is to beg the question. Numbers are used to count things, and the notion is fully expressed
by numerals.

To argue by induction that variable-free primitive-recursive terms reduce to numerals is to justify
finitary reasoning by an infinitary argument.
Not only is induction as a general axiom schema lacking any justification
other than an appeal to $\N$ as a completed infinity, but its
application to specific variable-free primitive-recursive terms
lacks a cogent justification.

We shall exhibit a primitive recursion and prove that it
does not terminate, thus disproving Church's Thesis from below and demonstrating
that finitism is untenable.

\bigskip

\hfill\vbox{\baselineskip=10pt\halign{\ninerm#\cr
The most extreme view, held by at least one mathematician\cr
at a respectable university, is that eventually a contradiction\cr
will be found even in elementary number theory.\hfil\cr}}\smallskip

\hfill\hbox{--- \ninerm P{\sevenrm AUL} \ninerm J.~C{\sevenrm OHEN}}

\vfill\eject

\chap2.{Logic}
This chapter explains the kind of proof we shall use, but the main
purpose is to discuss informally and qualitatively some of
the devices that later will be established formally in a theory.
Much of proof theory consists in algorithms for constructing or transforming proofs.
For reasons that will appear later, polynomially bounded algorithms are unproblematic,
but we formulate results involving exponential or superexponential bounds
as conditional results with the proviso that the algorithms terminate.
The exposition is based on~[Sh] but with greater use of special constants.

{\everymath={\rm} \everydisplay={\rm}

\section//   3.{Languages}
We use ``language'' to mean what is also called a first-order language.
The symbols of a language are \ii variables., \ii logical operators., \ii function
symbols., and \ii predicate symbols.. We use~s as a syntactical variable [Sh~\S1.3]
for symbols. Each symbol~s has a number~$I(s)$, called its \ii index..
Symbols of index 0, 1, 2, 3, and so forth, are called 0-ary, unary, binary,
ternary, and so forth. (For this reason, many writers use the neologism
``arity'' for the index.) A concatenation $s_1s_2\ldots s_\nu$
of symbols is an \ii expression., and $\nu$ is the \ii length.\/ of the
expression. The case $\nu=0$, the empty expression, is allowed.
We use u v as syntactical variables for expressions.

The \ii logical operators.\/ are $\Neg$ (\ii not.\/), $\vee$ (\ii or.\/),
and \hbox{$\exists$~(\ii there exists.\/)},
with $I(\Neg)=1$, $I({\vee})=2$, and $I(\exists)=2$.
We use x y z w as syntactical variables for variables, with $I(x)=0$.

The symbol~= (\ii equality.\/) is a binary predicate symbol.
It is important to distinguish
the use of~= as a symbol from its usual use in mathematics.

\goodbreak

The symbols already described are
common to all languages. A language is specified by giving in addition its \ii nonlogical symbols.,
consisting of function symbols (other than special constants, which will be discussed later)
and some predicate symbols other than~=, each with its index.
We shall consider only languages with finitely many nonlogical symbols, and for convenience we
require that the constant~0 is one of them.
We formally define a \ii language.\/ to be a concatenation of nonlogical symbols including~0, and
use~L as a syntactical variable for languages.
We say that u~is an expression \ii of.\/~L in case each nonlogical symbol occurring in~u occurs in~L.
We use f g as syntactical variables for function symbols, and p q for predicate symbols.
A \ii constant.\/ is a function symbol of index~0 and a \ii proposition letter.\/
(also called a \ii Boolean variable.\/)
is a predicate symbol of index~0. We use~e as a syntactical variable for constants.

\ii Terms.\/ are defined recursively as follows:
x is a term;
if $I(f)=\iota$ and $u_1$, \dots, $u_\iota$ are terms, so is
$fu_1\ldots u_\iota$.

Note that constants are terms. We use a b c d as syntactical variables for terms.

An \ii atomic formula.\/ is
$pa_1\ldots a_\iota$
where $I(p)=\iota$. \ii Formulas.\/ are defined recursively as follows:
an atomic formula is a formula;
if u and v are formulas, so are $\Neg u$, ${\vee}uv$, and $\exists xu$.

We use A B C D as syntactical variables for formulas.

Here is a summary of our syntactical variables to date:

\medskip

\halign{\qquad\qquad#\hfil&\quad for\quad#\hfil\cr
s&symbols\cr
u v &expressions\cr
x y z w&variables\cr
f g&function symbols\cr
e&constants\cr
a b c d&terms\cr
p q&predicate symbols\cr
A B C D&formulas\cr}

\medskip

We have been using prefix (or ``Polish'') notation, which makes general
discussions easier. But for readability we often use infix notation for binary symbols,
as in $A \orR B$ for ${\vee}AB$ and $a = b$ for ${=}ab$.
When we do this, groupers are necessary in general. We use parentheses to
group terms, brackets and braces to group formulas.
Infix symbols are associated from right to left.
We use $a\ne b$ to
abbreviate $\neg a = b$.
We call $\Neg A$ the \ii negation.\/ of~A and $A\orR B$ the \ii disjunction.\/
of A and~B.

The logical operators
\& (\ii and.\/), $\impP$ (\ii implies.\/), $\iffF$ (\ii if and only if.\/), and
$\forall$ (\ii for all.\/) are defined symbols as in [Sh~\S1.2];
we use the abbreviations

\medskip

\halign{\quad\qquad$#$\hfil&\quad for \quad $#$\hfil\cr
A\andD B&\Neg[\Neg A\orR \Neg B]\cr
A\impP B&\Neg A\orR B\cr
A\iffF B&[\Neg A\orR B]\andD[A\orR\Neg B]\cr
\forall x A&\Neg\exists x\Neg A\cr}

\medbreak

To avoid writing too many brackets, we make the convention that $\Neg$~and~$\exists x$
bind tightly and that ${\vee}$ and~\& bind more tightly than ${\impP}$ and~${\iffF}$.
For example, the formula
$\Neg\exists x A \impP \Neg B \orR C$
is parsed as
$\{ \Neg [ \exists xA ] \} \impP \{[ \Neg B ] \orR C \}$.

A \ii quantifier.\/ (\ii existential.\/ or \ii universal.\/ respectively) is
$\exists x$ or~$\forall x$.

We say that B is a \ii subformula.\/ of A in case B occurs in A; this terminology
differs from that of some writers.

An occurrence of x in A is \ii bound.\/ [Sh~\S2.4] in case it is
in a subformula $\exists xB$ of~A.
Otherwise, an occurrence of~x in~A is \ii free..
We say that b~is \ii substitutable for~$x$ in~$A$. in case no variable~y
occurring in~b becomes bound when each free occurrence of~x in~A is replaced by~b.
An expression is \ii variable free.\/ in case no variable occurs in it. If b~is
variable free then it it substitutable for~x in~A.
We use
$b_{x_1\ldots x_\nu}(a_1\ldots a_\nu)$
for the term obtained by replacing, for all~$\mu$ with $1\le\mu\le\nu$,
each occurrence of~$x_\mu$ in~b by~$a_\mu$.
We use
$A_{x_1\ldots x_\nu}(a_1\ldots a_\nu)$
for the formula, called an \ii instance.\/ of A,
obtained by replacing, for all~$\mu$ with $1\le\mu\le\nu$,
each free occurrence of~$x_\mu$ in~A by~$a_\mu$; when we do this, we always
assume that $a_\mu$~is substitutable for~$x_\mu$ in~A.
As an example of what can go wrong if this is not observed,
suppose that we replace
$\mit y$ by~$\mit x$ in $\exists\mit x[\mit x\ne\mit y]$; then we obtain
$\exists\mit x[\mit x\ne\mit x]$.

A formula is \ii open.\/ in case $\exists$ does not occur in it; it is \ii closed.\/
in case no variable occurs free in it.
Note that a formula is both open and closed if and only if
it is variable free. The \ii closure.\/ of~A, denoted by~$\bar A$, is
$\forall x_1\ldots\forall x_\iota A$
where $x_1$, \dots, $x_\iota$ are the variables occurring free in A
in the order of free occurrence.

The \ii index.\/ of a formula A, denoted by I(A), is the number of
distinct free variables occurring in it. If $I(A)=\iota$, A is called $\iota$-ary
(closed, unary, binary, etc.).
We use~$\phi$ as a syntactical variable for unary
formulas with free variable~$\mit x$, and we write $\phi(a)$ for~$\phi_{\mit x}(a)$.
(Requiring them all to have the same free variable ensures that $\phi_1\andD\phi_2$ is
unary.)
We also use~$\phi$ as a syntactical variable for unary predicate symbols; context will
make it clear which is intended.

\bigskip

\hfill\vbox{\baselineskip=10pt\halign{\ninerm#\hfil\cr
Die Mathematiker sind eine Art Franzosen:\cr
Redet man zu ihnen, so \umlaut ubersetzen sie es in ihre Sprache,\cr
und dann ist es alsobald ganz etwas anders.\cr}}\smallskip

\hfill\hbox{--- \ninerm J{\sevenrm OHANN} \ninerm W{\sevenrm OLFGANG VON} \ninerm G{\sevenrm OETHE}}

\section//   4.{Structures}
Before exploring a strange city, it is useful to look at a map. The map can give
some idea of where to go and what one might expect to find. But the city is alive
and growing, full of unpredictable events, while the map is dead and static.
The wise tourist does not identify the map with the city.

Semantics aids understanding if it is not taken literally: a structure gives
metaphoric meaning to a language. Let L be a language.
Then a \ii structure.\/ $\sigma$ for~L [Sh~\S2.5] consists of a non-empty set~U,
the \ii universe.\/ of~$\sigma$, and for each $\iota$-ary f in~L a function

\noj $\sigma(f):\; U^\iota \impP U$

\noin and for each $\iota$-ary p in~L a subset

\noj $\sigma(p)\subseteq U^\iota$

\noin Thus the meaning, in a structure, of a function
symbol is a function, and the meaning, in a structure, of a predicate symbol is
a predicate (i.e., a set). The elements of the universe~U are called \ii individuals..

The functions and predicates of a structure corresponding to the function symbols
and predicate symbols of~L are completely arbitrary. Recall that
=~does not occur in~L. The predicate assigned to it is not arbitrary; it is
the diagonal of~$U^2$. We define

\noj $\sigma({=}) = \{ (\xi_1,\xi_2) \in U^2 \;:\; \xi_1 = \xi_2 \, \}$

\noin Note that the first
occurrence of~= here is the equality symbol and the other occurrences are
the mathematical use of =.

The meaning, in a structure, of a variable is that it denotes some individual;
which individual is arbitrary, depending on which individual is assigned to it.
Let X~be the set of all variables.
An \ii assignment.\/ is a mapping $\alpha\colon X \impP U$.
The set of all assignments is~$U^X$.
The meaning, in a structure, of a term will depend on the
individuals assigned to the variables in it. We define $\sigma(b)$ as a mapping from
assignments~$\alpha$
to individuals; it is the individual denoted by~b in the structure~$\sigma$
when the variables in~b take the values assigned to them by~$\alpha$.
The definition is recursive:

\noj $\big(\sigma(x)\big)(\alpha)=\alpha(x)$

\noj $\big(\sigma(fa_1\ldots a_\iota)\big)(\alpha)=\big(\sigma(f)\big)
\Big(\big(\sigma(a_1)\big)(\alpha),\ldots,\big(\sigma(a_\iota)\big)(\alpha)\Big)$

\goodbreak

\noin where $I(f)=\iota$.

The meaning, in a structure, of a formula is that it is true or false, but the
truth value will depend on the individuals assigned
to the variables occurring free in it.
We define $\sigma(A)$ to be the set of assignments satisfying it, beginning with
atomic formulas:

\noj $\sigma(pa_1\ldots a_\iota)= \Big\lbrace\alpha \in U^X \,:\,
\Big(\big(\sigma(a_1)\big)(\alpha),\ldots,\big(\sigma(a_\iota)\big)(\alpha)\Big)
\in \sigma(p)\,\Big\rbrace$

\noin where $I(p)=\iota$.

If $\alpha$ is an assignment and $\xi$ is an individual, we let
$\alpha_{x,\xi}(x)=\xi$ and $\alpha_{x,\xi}(y)=\alpha(y)$
whenever y in a variable other than~x.
This alters $\alpha$ so that its value at~x is $\xi$ and leaves the assignment
of other variables unaltered.
Now let W~be a set of assignments. We define its \ii projection.\/ by

\noj $\pi_{x}W=\{\alpha\in U^X\;:\;\hbox{ there exists $\xi$ in U
such that $\alpha_{x,\xi}\in W$}\,\}$

\noin This consists of all assignments that become elements of~W just by changing
the assignment on~x. Also, the \ii complement.\/ of W is

\noj $\co W=U^X\setminus W$

\noin consisting of all assignments that are not in W.

Now we can recursively complete the definition of $\sigma(A)$:

\noj $\sigma(\Neg A)=\co{\sigma(A)}$

\noj $\sigma(A\orR B) = \sigma(A) \cup \sigma(B)$

\noj $\sigma(\exists xA) = \pi_x \sigma(A)$

\medskip

We say that an assignment~$\alpha$ \ii satisfies.\/ A in case $\alpha\in\sigma(A)$.
This depends only on the values that $\alpha$ assigns to the free variables of~A.
Similarly, if $\alpha$ and $\alpha'$ are two assignments that agree on all
the variables occurring in~b, then
$\big(\sigma(b)\big)(\alpha)=\big(\sigma(b)\big)(\alpha')$.

Let $x_1$, \dots, $x_\iota$ be, in some order, the variables occurring free
in~A. Then we define $\sigma_{x_1,\ldots,x_\iota}(A)$
to be the set of all $(\xi_1,\ldots,\xi_\iota)$ in~$U^\iota$ such that for all
assignments~$\alpha$ for which
$\alpha(x_1)=\xi_1$, \dots,~$\alpha(x_\iota)=\xi_\iota$
we have $\alpha\in\sigma(A)$.

Let \T\ be the set of all assignments (for a given universe U) and let
\F\ be the empty set of assignments.
If A is closed, either $\sigma(A)=\T$,
in which case we say that A~is \ii true.\/ in~$\sigma$, or $\sigma(A)=\F$,
in which case we say that A~is \ii false.\/ in~$\sigma$. We say that A~is
\ii valid.\/ in~$\sigma$ in case $\sigma(A)=\T$;
this is equivalent to saying that the closure of~A is true in~$\sigma$.
Two formulas A and~B are \ii equivalent in~$\sigma$. in case $A\iffF B$ is
valid in~$\sigma$, which is to say that they are satisfied by the same assignments;
i.e., $\sigma(A)=\sigma(B)$.

We say that $L'$ is an \ii extension.\/ of L in case every symbol of L is a
symbol of~$L'$. Let $L'$ be an extension of~L and let $\sigma$~be a structure
for~L. Then a structure~$\sigma'$ for~$L'$ is an \ii expansion.\/ of~$\sigma$
in case it has the same universe and $\sigma'(s)=\sigma(s)$ for all s in~L.

\goodbreak

Needless to say, structures will be used only as illustrations, not as part of the
formal development. Semantics is a topic in \ZFC, not an independent body of truths.

\bigskip

\hfill\vbox{\baselineskip=10pt\halign{\ninerm#\hfil\cr
He had bought a large map representing the sea\cr
\quad Without the least vestige of land:\cr
And the crew were much pleased when they found it to be\cr
\quad A map they could all understand.\cr
\cr
``What's the good of Mercator's North Poles and Equators,\cr
\quad Tropics, Zones, and Meridian Lines?''\cr
So the Bellman would cry: and the crew would reply\cr
\quad ``They are merely conventional signs!\cr
\cr
``Other maps are such shapes, with their islands and capes!\cr
\quad But we've got our brave Captain to thank''\cr
(So the crew would protest) ``that he's bought us the best---\cr
\quad A perfect and absolute blank!''\cr}}\smallskip

\hfill\hbox{--- \ninerm L{\sevenrm EWIS} \ninerm C{\sevenrm ARROLL}}

\section//   5.{Tautologies}
Now we study the propositional calculus---which ignores variables, equality,
function symbols, predicates symbols, and quantifiers---involv\-ing
only $\Neg$ and~$\vee$.
A formula is \ii elementary.\/ in case it is atomic or begins with~$\exists$.
Every formula is built up from elementary formulas by repeatedly forming negations
and disjunctions.
A \ii truth valuation.\/~V [Sh~\S3.1] is a function from the elementary formulas
to $\lbrace\T,\F\rbrace$.
We extend it to all formulas, still taking values either \T\ or~\F\
and retaining the notation~V, recursively
as follows:

\noj $V(\Neg A) = \T$\qquad\hskip3.85pt if and only if\quad $V(A)=\F$

\noj $V(A \orR B) =\T$\quad if and only if\quad $V(A)=\T$\quad
or\quad $V(B)=\T$

\noin Then A is a \ii tautology.\/ in case $V(A)=\T$ for all truth valuations,
and A is a \ii tautological consequence.\/ of $A_1$, \dots,~$A_\nu$
in case $A_1\andD\cdots\andD A_\nu\impP A$
is a tautology, which is equivalent
to saying that for all truth valuations~V, if $V(A_\mu)=\T$
for all~$\mu$ with $1\le\mu\le\nu$, then $V(A)=\T$.
This is a purely syntactical notion, not involving the semantics of structures.
But note that if \hbox{$A_1$, \dots,~$A_\nu$}
are true in a structure~$\sigma$ and A~is
a tautological consequence of them, then A~is true in~$\sigma$,
since on closed formulas $A\mapsto\sigma(A)$ is a truth valuation.

The formulas A and B are \ii tautologically equivalent.\/ in case $A\iffF B$
is a tautology.

To check whether A~is a tautology it suffices to consider
truth valuations defined just on the subformulas of~A.
And one might as well replace
all the elementary formulas in it by Boolean variables, with different elementary
formulas being replaced by different Boolean variables. If there are $\nu$
Boolean variables then there are $2^\nu$ relevant truth valuations, since each
Boolean variable independently may be assigned the truth value \T\ or~\F.
For each truth valuation~V one quickly evaluates it on the formula by using the
two rules above; this is called the method of truth tables. If $\nu$ is
large, $2^\nu$ is enormous and the method is infeasible.

\section//   6.{Equality}
The next level up from the propositional calculus is what may be called the
functional calculus, in which variables, equality,
predicate symbols, and function symbols play a role but $\exists$~is still ignored.

An \ii identity axiom.\/ is

\noj $x=x$

\noin and an \ii equality axiom.\/ is

\noj $x_1 = y_1 \and \cdots \and x_\iota = y_\iota \and px_1\ldots x_\iota
\imp py_1\ldots y_\iota$

\noin or

\noj $x_1 = y_1 \and \cdots \and x_\iota = y_\iota \imp fx_1\ldots x_\iota
= fy_1\ldots y_\iota$

\noin where p and f have index $\iota$. Identity and equality axioms are valid in
every structure.

A formula is a \ii quasitautology.\/ in case it is a tautological consequence
of instances of identity and equality axioms, and A~is a \ii quasitautological
consequence of.\/ $A_1$, \dots,~$A_\nu$ in case $A_1\andD\cdots\andD A_\nu
\impP A$ is a quasitautology.

\goodbreak

\Th//  1. {\it The formulas

\noj $a=a$

\noj $a=b \imp b=a$

\noj $a=b \and b=c \imp a=c$

\noin are quasitautologies.}

\pf The first is an instance of an identity axiom, and the others are
tautological consequence of the instances

\noj $a=b \and a=a \and a=a \imp b=a$

\noj $a=a \and b=c \and a=b \imp a=c$

\noin of the following equality axiom for equality

\noj $x_1=y_1\and x_2=y_2 \and x_1=x_2 \imp y_1=y_2$

\noin and of instances of identity axioms. \bul

\bigskip

\vbox{
\hfill\vbox{\baselineskip=10pt\halign{\ninerm#\hfil\cr
Things equal to the same thing are equal.\cr}}\smallskip

\hfill\hbox{--- \ninerm E{\sevenrm UCLID}}

}

\section//   7.{Special constants}
Now we come to the full predicate calculus, quantifiers and all.
A formula of the form $\exists xA$ is an \ii instantiation., and
$\forall xA$ is a \ii generalization..

To each closed instantiation $\exists xA$ we associate a constant
$c^{\phantom0}_{\exists xA}$, called the \ii special constant for.\/~$\exists xA$,
and $\exists xA$ is its \ii subscript.. We use~r as a syntactical
variable for special constants. If r~is the special constant for~$\exists xA$,
its \ii special axiom.\/ is

\noj $\exists xA\impP A_x(r)$

\noin We say that u~is \ii plain.\/ in case no special constant occurs in it.

The \ii level.\/ [Sh~\S4.2] of a special constant is defined recursively as follows.
If A~is plain, the level of $c^{\phantom0}_{\exists xA}$ is~1; otherwise its level is
1~plus the maximal level of the special constants occurring in~A.

Mathematicians use special constants all the time with no special mention.
For example, someone studying odd perfect numbers in the hope of showing that they
do not exist may call such a number~$\mit n$ and
proceed to derive properties of~$\mit n$. Then~``$\mit n$'' is a special constant.
The alternative would be to enclose the entire discussion in a huge formula of the form
\hbox{$\exists\mit n[\mit n \' is an odd perfect number' \impP \cdots]$}.
The practice is to use the same notation for
special constants as for variables, bearing in mind that the special constants are
fixed throughout the discussion.

We regard special constants as logical symbols; they are a device to facilitate
the handling of quantified formulas in proofs. The \ii nonlogical symbols.\/
are the function symbols other than special constants and the predicate symbols
other than equality; the other symbols (variables, logical operators, equality,
and special constants) are the \ii logical symbols..

Semantically, let $\sigma$~be
a structure for L. Well-order the universe~U of~$\sigma$.
We expand~$\sigma$ to be a structure $\hat\sigma$ for the language schema
obtained by adjoining all special constants, for formulas all of whose nonlogical
symbols are in~L, by recursion on the level, as follows.
Let r~be the special constant for
$\exists xA$ and suppose that $\hat\sigma(r')$ has been defined for all~$r'$
of strictly lower level.
If there exists an individual~$\xi$ and an assignment~$\alpha$
with $\alpha(x)=\xi$ such that $\alpha\in\hat\sigma(A)$,
let $\hat\sigma(r)$ be the least such individual;
otherwise, let $\hat\sigma(r)$ be~$\sigma(0)$.
Then the special axioms are true in~$\hat\sigma$.

We define the notion that v \ii appears in.\/ u recursively as follows:
if v~occurs in~u, then v~appears in~u; if v occurs in the subscript of
a special constant that appears in~u, then v~appears in~u.
(That is, v~appears in~u if and only if it is there when all subscripts of special
constants are written out; ``appears in'' is the transitive closure of
``occurs in or occurs in the subscript of a special constant occurring in''.)

\bigskip

\hfill\vbox{\baselineskip=10pt\halign{\ninerm#\hfil\cr
Generalizations are seldom or ever true,\cr
and are usually utterly inaccurate.\cr}}\smallskip

\hfill\hbox{--- \ninerm A{\sevenrm GATHA} \ninerm C{\sevenrm HRISTIE}}

\section//   8.{Proofs}
Lemma 1 of [Sh~\S4.3] is the key to a notion of proof that is fully formal and yet
close to rigorous arguments in mathematical practice.
This section is devoted to an explication of this notion.

We shall say that a \ii theory.\/ is a concatenation of formulas, and use~T
as a syntactical variable for theories. The formulas occurring
in~T are the \ii nonlogical axioms.\/ of~T.
(One of the virtues of prefix notation is that an expression
can be a concatenation of formulas in only one way, so we can define a theory
to be an expression rather than a finite sequence of formulas.)
In other words, all our theories are finitely axiomatized. If infinitely many nonlogical
axioms are specified in some way, we call this a \ii theory schema..
The \ii language of.\/~T, denoted by~L(T), is the concatenation of the nonlogical
symbols occurring in~T, together with~0.
A theory is \ii open.\/ in case all its nonlogical axioms are open and plain. By $T[A]$
is meant the theory obtained from~T by adjoining~A as a new nonlogical axiom.

Semantically, a \ii model.\/ of T is a structure for L(T) in which all of the nonlogical
axioms are valid.

A \ii substitution axiom.\/ is $A_x(a)\impP\exists xA$. (Substitution axioms are
valid in every structure.)
By a \ii substitution formula., \ii identity formula., or \ii equality formula.\/
we mean respectively a closed instance of a substitution axiom, identity axiom, or
equality axiom.

Given a theory T, we say that a formula is \ii in $\Delta(T)$. in case it is
a formula of one of the following five forms:

\nok// 24. special axiom:
$\exists xA\impP A_x(r)$,\quad r the special constant for $\exists xA$

\nok// 25. substitution formula: $A_x(a)\impP\exists xA$

\nok// 26. identity formula

\nok// 27. equality formula

\nok// 28. closed instance of a nonlogical axiom of T

\medskip

We call $(/ 24)$--$(/ 27)$ \ii logical axioms..

Shoenfield's Lemma 1 asserts that if there is a proof of A in~T and $A'$ is a closed
instance of~A, then $A'$~is a tautological consequence of formulas
in $\Delta(T)$ (and conversely).
In this lemma, ``proof'' refers to the notion as formulated
in [Sh~\S2.6], but we shall take this result as a definition.

Use $\pi$ as a syntactical variable for concatenations of formulas, and call~A a \ii formula of.\/~$\pi$
in case $\pi$ is of the form $\pi_1A\pi_2$.
Call $\pi$ a \ii proof.\/ in~T in case each of its formulas is in~$\Delta(T)$ or is
a tautological consequence of strictly preceding formulas.
Every formula in a proof is closed.
A proof~$\pi$ is a \ii proof of.\/~A, denoted by $\pi\vdash_T A$,
in case the closure of~A
is a formula of the proof, in which case A~is a \ii theorem.\/ of~T, denoted by
$\vdash_T A$.

There is no requirement that a theorem of~T be plain or a formula
of~$L(T)$. But extraneous symbols can be eliminated from proofs.
We call a mapping $A\mapsto A^+$ defined on formulas of~L a \ii homomorphism on.\/ L in
case $[\Neg A]^+$ is $\Neg A^+$ and $[A\orR B]^+$ is $A^+\orR B^+$
and $[\exists xA]^+$ is $\exists xA^+$. Thus it suffices to give the
values of a homomorphism on atomic formulas.
Notice that homomorphisms preserve tautological consequence.

\goodbreak

\Th//  2. {\it If\/ $\pi\vdash_T A$, then there
is an expression\/~$\pi^+$ of\/~$L(T[A])$ such that\/ $\pi^+\vdash_T~A$.}

\pf Let $u^+$ be obtained by replacing, everywhere it appears,
each atomic formula beginning with a predicate
symbol not in~$L(T[A])$
by~$0=0$ and then replacing each term beginning with a function symbol that is not
in~$L(T[A])$ by~0. Then $u^+$ is an expression of~$L(T[A])$ and
$B\mapsto B^+$ is a homomorphism. If B is a special axiom, substitution axiom, identity
formula, or closed instance of a nonlogical axiom of~T, so is~$B^+$. If B~is an
equality formula for a symbol of~$L(T[A])$, so is~$B^+$, and otherwise it is a
tautological consequence of the identity formula~$0=0$. \bul

\medskip

This is a pattern of proof that will often be employed. We construct $u\mapsto u^+$
such that $B\mapsto B^+$ is a homomorphism, and so preserves
tautological consequence. Then we need only examine formulas $(/ 24)\hbox{-}(/ 28)$
and verify whatever properties are stated in the theorem under consideration.

A theory $T_1$ is an \ii extension.\/ of the theory T in case
every nonlogical axiom of~$T_1$ is a
theorem of~T. It is a \ii conservative extension.\/ in case every formula of~$L(T)$
that is a theorem of~$T_1$ is a theorem of~T.
A theory~T is \ii consistent.\/ in case $0\ne0$ is not a theorem of~T.
Clearly, if $T_1$~is a conservative extension of~T and T~is consistent,
so is~$T_1$.

\Th//  3. {\it If\/ $ \vdash_T A$ then\/ $T[A]$ is a conservative extension of\/~\ro T.}

\pf Let the closure of A be $\forall x_1\ldots\forall x_\iota A$. By hypothesis,
this is a formula of a proof in~T.
Since the only new formulas in $\Delta(T[A])$ are closed instances of~A,
we need to show that such a closed instance~B is a formula of a proof in~T.
Now a closed formula of the form $\forall xC \impP C_x(a)$
is a tautological consequence of the substitution formula
\hbox{$\Neg C_x(a)\impP\exists x\Neg C$},
and so is a logical theorem. Apply this remark $\iota$~times.~\bul

\medskip

The rather trivial extension of T of the form $T[A]$, where A is a plain formula and a
theorem of~T, is called a \ii t-extension of.\/~T.

\bigskip

\hfill\vbox{\baselineskip=10pt\halign{\ninerm#\hfil\cr
There's nothing you can't prove\cr
if your outlook is only sufficiently limited.\cr}}\smallskip

\hfill\hbox{--- \ninerm D{\sevenrm OROTHY} L. S{\sevenrm AYERS}}

\goodbreak

\section//   9.{Some logical theorems}
We say that A~is a \ii logical theorem., denoted by $\vdash A$,
in case it is a theorem of
the theory with no nonlogical axioms, and A and~B are \ii logically equivalent.\/ in case
$\vdash A\iffF B$.

Let r be the special constant for $\exists xA$. Then $\exists xA\impP A_x(r)$
is a special axiom and the formula $A_x(r)\impP\exists xA$ is a substitution formula,
so $\exists xA$ and $A_x(r)$ are logically equivalent. By repeated use of this,
we see that any closed formula is logically equivalent to a variable-free formula
(but one that may contain special constants).
The use of special constants reduces reasoning in the predicate calculus to
reasoning in the functional calculus.

The following is the \ii equality theorem.\/ [Sh~\S3.4].

\Th//  4. {\it A closed formula of the form}

\nok// 29. $ a_1=b_1 \and \cdots \and a_\iota=b_\iota \and
A_{x_1\ldots x_\iota}(a_1\ldots a_\iota) \imp
A_{x_1\ldots x_\iota}(b_1\ldots b_\iota)$

\noin{\it is a logical theorem.}

\pf The proof is by induction on the \ii height.\/ of~A (the number of occurrences
of logical operators in~A). If A~is atomic, $(/ 29)$ is an equality formula.
If A~is $\Neg B$, the result holds by the induction hypothesis, using
Theorem/  1 to reverse the role of the a's and b's. If A~is $ B\orR C$,
the result holds by the induction hypothesis. If A is~$\exists xB$, let r~be
the special constant for $\exists xB_{x_1\ldots x_\iota}(a_1\ldots a_\iota)$.
By the induction hypothesis we have

\medskip

\centerline{$ a_1=b_1 \and \cdots \and a_\iota=b_\iota \and
B_{xx_1\ldots x_\iota}(ra_1\ldots a_\iota) \imp
B_{xx_1\ldots x_\iota}(rb_1\ldots b_\iota)$.}

\noin But $(/ 29)$ is a tautological consequence of this, the special
axiom for~r, and the substitution formula
$ B_{xx_1\ldots x_\iota}(rb_1\ldots b_\iota)\impP \exists x
B_{x_1\ldots x_\iota}(b_1\ldots b_\iota)$.~\bul

\medskip

Let the closure $\bar A$ of A be $\forall x_1\ldots\forall x_\iota A$.
Let $A^1$ be~$\bar A$ and, recursively for $1\le\kappa<\iota$ let $A^{\kappa+1}$
be $A^\kappa_{x_\kappa}(r_\kappa)$ where $r_\kappa$ is the special constant
for~$\exists x_\kappa \Neg A^\kappa$. Denote $A^\iota$ by~$A^\circ$, and call this
A~\ii with frozen variables..

\Th//  5. {\it $\bar A\iffF A^\circ$ is a logical theorem. In any theory\/ \ro T,
$\vdash_TA$ if and only if\/ $\vdash_TA^\circ$.}

\pf $\vdash\exists x_\kappa\Neg A^\kappa\iffF\Neg A^\kappa(r_\kappa)$ since
the forward direction is a special axiom and the backward direction is a
substitution formula, but $\forall x_\kappa A^\kappa\iffF A^{\kappa+1}$ is a tautological
consequence of this. The second statement holds since, by definition, A~is a
theorem of~T if and only if $\bar A$~is.\bul

\medbreak
\goodbreak

The following is the \ii equivalence theorem.\/ [Sh~\S3.4].

\Th//  6. {\it Let\/ $ A'$ be obtained from\/ \ro A by replacing an occurrence
of a subformula\/~\ro B of\/~\ro A by\/~$ B'$. If\/ $\vdash_T B\iffF B'$
then\/ $\vdash_T A\iffF A'$.}

\pf By induction on the height of A, it suffices to consider the cases that A is
$\Neg B$, $ B\orR C$, $ C\orR B$, or $\exists xB$. Freezing variables,
we assume that these formulas are closed. The first three
cases hold by tautological equivalence, so suppose that A is $\exists x B$,
$ A'$ is $\exists xB'$, and $\vdash_T B\impP B'$. From~A we have
$ B_x(r)$, where r is the special constant for $\exists xB$, and so $ B'_x(r)$ and hence
$\exists xB'$; that is, $\vdash_T A\impP A'$. The other direction is the same.~\bul

\medskip

Next is the \ii deduction theorem.\/ [Sh~\S3.3].

\Th//  7. {\it Let\/ \ro C be closed. If\/ $\vdash_{T[C]}A$ then\/ $\vdash_T C\impP A$.}

\pf If B is a special axiom, substitution formula, identity formula, equality formula,
or nonlogical axiom of T[C], then $\vdash_T C\impP B$. If B is a tautological
consequence of \hbox{$B_1$, \dots,~$B_\nu$}, then $C\impP B$ is a tautological consequence
of $C\impP B_1$, \dots,~$C\impP B_\nu$. \bul

\medskip

The negation form of a formula is obtained by pushing each $\Neg$ repeatedly to the
right using the rules of the predicate calculus. The precise definition is complicated
by the fact that we chose, following [Sh~\S2.4], to regard~$\forall$ and~$\&$ as
defined symbols. This is convenient in many other contexts but a nuisance here.
Recall that $\forall x A$ abbreviates $\Neg\exists x\Neg A$ and that $ A\andD B$
abbreviates $\Neg{\vee}\Neg A\Neg B$. Call an
occurrence of~$\Neg$ \ii proper.\/ in case it is not displayed in one of these two forms.
Now repeatedly
make the following replacements, where the displayed occurrences of~$\Neg$ are proper.

\medskip

\halign{\qquad\qquad$#$\hfil&\qquad$#$\hfil\cr
replace:&by:\cr
\noalign{\vskip5pt}
\Neg\exists xA&\forall x\Neg A\cr
\Neg\forall xA&\exists x\Neg A\cr
\Neg[A\hskip1pt\orR\hskip1pt B]&\Neg A\andD\Neg B\cr
\Neg[A\andD B]&\Neg A\hskip1pt\orR\hskip1.4pt\Neg B\cr
\Neg\Neg A&A\cr
}

\medskip

\noin The result of making all these replacements in the subformulas of
a formula is the \ii negation form.\/ of the formula.
If A~is in negation form, then every proper occurrence of~$\Neg$ immediately precedes an
atomic formula.

\Th//  8. {\it Let\/ $ A'$ be the negation form of\/ \ro A. Then\/
$\vdash A\iffF A'$.}

\pf We have $\vdash A\iffF\Neg\Neg A$, so
$\vdash \Neg\exists xA\iffF\Neg\exists x\Neg\Neg A$ by the equivalence theorem,
but $\Neg\exists x\Neg\Neg A$ is $\forall x\Neg A$. This proves the equivalence
of the first replacement in the definition of negation form. The
others are tautological equivalences, so the result holds by the
equivalence theorem.~\bul

\medskip

\Th//  9. {\it The following are logical theorems
provided that\/ $x$ is not free in\/~$C$.}

\nok// 30. $\exists x B\,\orR C\hskip.2pt \iff \exists x[B\,\orR C\hskip.2pt]$

\nok// 31. $C\,\orR \exists x B\hskip.2pt \iff \exists x[C\,\orR B\hskip.2pt]$

\nok// 32. $\exists x B\andD C \iff \exists x[B\andD C]$

\nok// 33. $C\andD \exists x B \iff \exists x[C\andD B]$

\nok// 34. $\forall x B\,\orR C\hskip.2pt \iff \forall x[B\,\orR C\hskip.2pt]$

\nok// 35. $C\,\orR \forall x B\hskip.2pt \iff \forall x[C\,\orR B\hskip.2pt]$

\nok// 36. $\forall x B\andD C \iff \forall x[B\andD C]$

\nok// 37. $C\andD \forall x B \iff \forall x[C\andD B]$

\pf Freeze the variables in~$(/ 30)$.
For the forward direction,
let r~be the special constant for~$\exists xB$. Then $B_x(r)\orR C$ is
a tautological consequence of $\exists x B\orR C$ and
the special axiom for~$\exists xB$, and
$B_x(r)\orR C\impP\exists x[B\orR C]$ is a substitution formula (since x~is
not free in~C). For the backward direction, let $r'$~be the special constant
for~$\exists x[B\orR C]$. Then \hbox{$B_x(r')\orR C$} is
a tautological consequence of $\exists x[B\orR C]$ and its special axiom
(again since x~is not free in~C),
and $B_x(r')\orR C\impP \exists x B\orR C$ is a tautological consequence
of the substitution formula $B_x(r')\impP \exists xB$. The other
seven cases are similar (for~$\forall xB$, use the special constant
for~$\exists x\Neg B$).\bul

\medskip

A \ii variant.\/ of A [Sh~\S3.4] is a formula obtained from A by replacing a subformula
$\exists xB$ by $\exists yB_x(y)$, where y~is not free in~B, zero or more times.
Variants are useful for avoiding colliding variables. We can find a variant of~A
such that no variable in~u occurs bound in~A; then each term occurring in~u
is substitutable for each free variable of~A.

Call a formula \ii adjusted.\/ in case different occurrences of quantifiers have
different bound variables and no variable occurs both free and bound. Every formula
has an adjusted variant.

The relation ``variant of'' is an equivalence relation. Semantically, variants are
equivalent to each other in any structure. The following is the \ii variant theorem.\/
[Sh~\S3.4].

\goodbreak

\Th//  10. {\it Let\/ $A'$ be a variant of\/~\ro A. Then $A\impP A'$ is a logical
theorem.}

\pf Freezing variables, we assume that A and $A'$ are closed.
The proof is by induction on the height of~A.

If A is atomic, then $A'$ is A, and $A\impP A$ is a tautology.

If A is $\Neg B$, so that $A'$ is $\Neg B'$ where B is a variant of $B'$, then we have
$B'\impP B$ by the induction hypothesis (since if B is a variant of~$B'$,
then $B'$~is a variant of~B), and $A\impP A'$ is a
tautological consequence of this.

If A is $B\orR C$, so that $A'$ is $B'\orR C'$ where $B'$ is a variant of B
and $C'$ is a variant of C, then $A\impP A'$ is a tautological consequence
of $B\impP B'$ and $C\impP C'$, which are logical
theorems by the induction hypothesis.

If A is $\exists xB$, so that $A'$ is $\exists x'B'$ where $B'$ is a variant of B,
let r be the special constant for A. Then $A\impP B_x(r)$ is a special axiom,
the formula $B_x(r)\impP B'_{x'}(r)$ is a logical
theorem by the induction hypothesis, and
$B'_{x'}(r)\impP A'$ is a substitution formula. But $A\impP A'$ is a
tautological consequence of these formulas. \bul

\medskip

A formula is in \ii prenex form.\/ in case it consists of a concatenation of
quantifiers (the \ii prefix.\/) followed by an open formula (the \ii matrix.\/).
Every formula is
logically equivalent to one in prenex form: take the negation form of an
adjusted variant and replace the left hand side of one $(/ 30)$--$(/ 37)$
by the right hand side, as often as possible.
The resulting formula is a \ii prenex form of.\/ the formula.
%The following definitions
%are made recursively: a formula is a \hbox{$\Sigma_0$-formula} in case it is open; it
%is a \hbox{$\Pi_0$-formula} in case it is open; it is a
%\hbox{$\Sigma_{\nu+1}$-formula} in case it
%consists of one or more existential quantifiers followed by a \hbox{$\Pi_\nu$-formula};
%it is a \hbox{$\Pi_{\nu+1}$-formula} in case it
%consists of one or more universal quantifiers followed by a \hbox{$\Sigma_\nu$-formula}.

The following are tautologies:

\noj $A\orR[B\andD C] \iff[A\orR B] \and [A\orR C]$

\noj $[B\andD C]\orR A \iff[B\orR A] \and [C\orR A]$

\goodbreak

\noin A formula is in \ii conjunctive form.\/ in case it contains no subformula
that is the left hand side of one of these equivalences. Unlike the previous
procedures,
repeatedly replacing the left hand side by the right hand side
may lead to an exponentially longer formula.

A formula is in \ii normal form.\/ in case it is closed and
in prenex, negation, and conjunctive form.

\section//   10.{Predicate extensions}
Mathematics would not get very far without the freedom to define new concepts.
Let D~be a plain formula of~T and let
the free variables of D, in the order of free occurrence, be $x_1$, \dots,~$x_\iota$.
Let $p^{\phantom0}_D$, the \ii predicate symbol for.\/~D, be
a predicate symbol of index~I(D). Then

\nok// 38. $ px_1\ldots x_\iota \iff D$

\noin where p is the predicate symbol for D, is the
\ii defining axiom of a predicate symbol.\/ [Sh~\S4.6]. The formula~D
is the \ii definiens.. Let $T_1$ be the
theory obtained from~T by adjoining~($/ 38$) as a new
nonlogical axiom; we call~$T_1$ a \ii p-extension of.\/~T.

Semantically, we can expand a model [Sh~\S2.6]
$\sigma$ of T to a model $\sigma_1$ of~$T_1$
by letting $\sigma_1(p)$ be~$\sigma_{x_1\ldots x_\iota}(D)$.

Define $u^+$ recursively, working from right to left,
by replacing each variable-free $pa_1\ldots a_\iota$
everywhere it appears in~u by
$D_{x_1\ldots x_\iota}(a_1^+\ldots a_\iota^+)$.
Then $A\mapsto A^+$ is a homomorphism.

\Th//  11. {\it Let\/ $T_1$ be a p-extension of\/ \ro T.
If\/ $\pi\vdash_{T_1}A$, then\/
$\pi^+\vdash_TA^+$.
The theory\/~$T_1$ is a conservative extension of\/~\ro T.}

\pf If B is a special axiom, substitution formula, identity formula,
or closed instance of a nonlogical axiom of~T, so is~$B^+$.

If B is a closed instance of the new axiom~$(/ 38)$, then $B^+$ is a tautology
of the form \hbox{$C\iffF C$}.

If B~is an equality formula for a function symbol or for a predicate
symbol other than the new symbol~p, so is~$B^+$.
If B~is the equality formula

\noj $a_1=b_1 \and \cdots \and a_\iota=b_\iota \and pa_1\ldots a_\iota \imp
pb_1\ldots b_\iota$

\noin then $B^+$ is

\noj $a_1^+=b_1^+ \and \cdots \and a_\iota^+=b_\iota^+ \and \hskip-1pt
D_{x_1\ldots x_\iota}(a_1^+\ldots a_\iota^+) \imp
D_{x_1\ldots x_\iota}(b_1^+\ldots b_\iota^+)$

\noin which is a logical theorem by the equality theorem.

The last statement of the theorem holds since
for a formula~A of~L(T), the formula~$A^+$ is A itself.~\bul

\medbreak

Let $L'$ extend L. A \ii translation from\/ $L'$ into.\/ L is a homomorphism
$A\mapsto A^*$ on~$L'$ such that each~$A^*$ is a formula of~L and the mapping is
the identity on formulas of~L.

Note that $A^+$ may not be a formula of L(T), even if A is plain: it may
contain non-closed atomic formulas beginning with~p, and p~may appear in subscripts
of special constants. Define a translation, unique up to variants,
as follows [Sh~\S4.6].
If A~is the formula $pa_1\ldots a_\iota$, choose a variant~$D'$ of~D in which
no variable of~A occurs bound and let $A^*$ be
$D'_{x_1\ldots x_\iota}(a_1\ldots a_\iota)$; for other atomic formulas~A, let
$A^*$ be~A.

\Th//  12. {\it Let\/ $T_1$ be a p-extension of\/ $T$. If\/ $A$ is a formula
of\/~$L(T_1)$, its translation~$A^*$ is a formula of\/~$L(T)$, and\/
$\vdash_{T_1} A\iffF A^*$. We have\/ $\vdash_{T_1}A$ if and only if\/
$\vdash_T A^*$.}

\pf The translation $A$ is by definition a formula of L(T). To prove that

\noj $\vdash_{T_1} A\iffF A^*$

\noin it suffices by the equivalence theorem to prove that

\nok// 39. $\vdash_{T_1} pa_1\ldots a_\iota
\iff D'_{x_1\ldots x_\iota}(a_1\ldots a_\iota)$

\noin Then

\nok// 40. $\vdash_{T_1} D_{x_1\ldots x_\iota}(a_1\ldots a_\iota)\iff
D'_{x_1\ldots x_\iota}(a_1\ldots a_\iota)$

\noin by the variant theorem, and $(/ 39)$ is a tautological consequence
of~$(/ 40)$ and a closed instance of the defining axiom~$(/ 38)$.

The last statement holds since $T_1$ is a conservative extension of T.\bul

\medbreak

\section//   11.{Function extensions}
Let $x_1$, \dots, $x_\iota$ be the free variables, in the order of free occurrence,
of~$\exists yD$, where D~is a plain formula of~T.
Let~$f^{\phantom0}_D$, the \ii function symbol for.\/~$\exists yD$, be
a function symbol of index~$I(\exists yD)$. The \ii existence condition.\/ is

\nok// EC. $\exists yD$

\noin and the \ii uniqueness condition.\/ is

\nok// UC. $D \and D_y(y') \imp y=y'$

\noin where $y'$ does not occur in D. If EC and UC are theorems of T, then

\nok// 41. $fx_1\ldots x_\iota = y \iff D$

\noin where f is the function symbol for $\exists yD$,
is the \ii defining axiom of a function symbol.\/ [Sh~\S4.6],
with \ii definiens.\/~D. Let $T_1$ be the
theory obtained from~T by adjoining~$(/ 41)$
as a new nonlogical axiom; we call~$T_1$ an \ii f-extension of.\/~T.

Semantically, we can expand a model $\sigma$ of T to a model $\sigma_1$ of~$T_1$
by letting $\sigma_1(f)$ be~$\sigma_{x_1\ldots x_\iota y}(D)$.
This is a function since EC and~UC are valid in~$\sigma$.

Define $u^+$ recursively, working from right to left,
by replacing each variable-free
$fa_1\ldots a_\iota$ everywhere it appears in~u by the special constant for
$\exists yD_{x_1\ldots x_\iota}(a_1^+\ldots a_\iota^+)$.
Then $A\mapsto A^+$ is a homomorphism.

\Th//  13. {\it Let\/ $T_1$ be an f-extension of\/ \ro T.
If\/ $\pi\vdash_{T_1}A$, then\/ $\pi^+\vdash_TA^+$.
The theory\/~$T_1$ is a conservative extension of\/~\ro T.}

\pf If B is a special axiom, substitution formula, identity formula,
or a closed instance of a nonlogical axiom of~T, so is $B^+$.

If~B is a closed instance of the new axiom~$(/ 41)$, and so of the form

\noj $ fa_1\ldots a_\iota=b \iff D_{x_1\ldots x_\iota y}(a_1\ldots a_\iota b)$

\noin then $B^+$ is

\nok// 42. $ r=b^+ \iff D_{x_1\ldots x_\iota y}(a_1^+\ldots a_\iota^+ b^+)$

\noin where r is the special constant for $\exists y
D_{x_1\ldots x_\iota}(a_1^+\ldots a_\iota^+)$.
By EC and the special axiom for~r,

\noj $\vdash_T D_{x_1\ldots x_\iota y}(a_1^+\ldots a_\iota^+r)$

\noin Hence the forward direction of $(/ 42)$ holds by the equality theorem,
and the backward direction by~UC.

If B~is an equality formula for a predicate symbol or for a function symbol other
than the new~f, so is~$B^+$. If B~is the equality formula

\noj $ a_1=b_1 \and \cdots \and a_\iota=b_\iota \imp fa_1\ldots a_\iota =
f\hskip.1pt b_1\ldots b_\iota$

\noin then $B^+$ is

\nok// 43. $ a_1^+=b_1^+ \and \cdots \and a_\iota^+=b_\iota^+ \imp r = r'$

\noin where r is the special constant for
$\exists yD_{x_1\ldots x_\iota}(a_1^+\ldots a_\iota^+)$
and $ r'$ is the special constant for
$\exists yD_{x_1\ldots x_\iota}(b_1^+\ldots b_\iota^+)$.
By EC and the special axiom for r,

\noj $\vdash_T \ D_{x_1\ldots x_\iota y}(a_1^+\ldots a_\iota^+ r)$

\noin Hence, by the equality theorem,

\nok// 44. $\vdash_T \ a_1^+=b_1^+ \and \cdots\and a_\iota^+=b_\iota^+
\imp D_{x_1\ldots x_\iota y}(b_1^+\ldots b_\iota^+ r)$

\noin By a closed instance of EC and the special axiom for $ r'$,

\nok// 45. $\vdash_T \ D_{x_1\ldots x_\iota y}(b_1^+\ldots b_\iota^+r')$

\noin But $(/ 43)$ is a tautological consequence of $(/ 44)$, $(/ 45)$, and the
closed instance

\noj $D_{x_1\ldots x_\iota y}(b_1^+\ldots b_\iota^+ r)\and
D_{x_1\ldots x_\iota y}(b_1^+\ldots b_\iota^+r') \imp r=r'$

\noin of UC.

The last statement of the theorem holds since
for a formula~A of~L(T), the formula~$A^+$ is A itself.~\bul

\medbreak
\goodbreak

Now define a translation $A\mapsto A^*$ as in [Sh~\S4.6]. Let A be an atomic formula
of~$L(T_1)$; we define~$A^*$ by recursion on the number of occurrences of~f in~A.
If there are none, let $A^*$ be~A. Otherwise, let B~be the atomic formula
obtained by replacing the rightmost occurrence of a term beginning with~f by~z,
where z~does not occur in~A, so that

\nok// 46. $A \quad is \quad B_z(fa_1\ldots a_\iota)$

\noin Then f does not occur in $a_1\ldots a_\iota$ and B has one less occurrence of f
than A, so $B^*$ is already defined by recursion. Let

\nok// 47. $A^* \quad be \quad
\exists z[D'_{x_1\ldots x_\iota y}(a_1\ldots a_\iota z)\and B^*]$

\noin If A~is a formula of $L(T_1)$, we call~$A^*$ the
\ii translation.\/ of~A.

\Th//  14. {\it Let\/ $T_1$ be an f-extension of\/ $T$. If\/ $A$ is a formula
of\/~$L(T_1)$, its translation~$A^*$ is a formula of\/~$L(T)$, and\/
$\vdash_{T_1} A\iffF A^*$. We have\/ $\vdash_{T_1}A$ if and only if\/
$\vdash_T A^*$.}

\pf That $A^*$ is a formula of L(T) follows from its definition, by
induction on the number of occurrences of~f in~A. By the equivalence theorem,
it suffices to prove that \hbox{$\vdash_{T_1} A\iffF A^*$} for A~an atomic formula
of~$L(T_1)$. Again, the proof is by induction on the number of occurrences of~f
in~A. If there are none, $A\iffF A^*$ is a tautology. Using the notation
introduced above, the induction hypothesis $(/ 46)$, and the equivalence theorem to
replace~$B^*$ by~B in~$(/ 47)$, we need to prove

\nok// 48. $B_z(fa_1\ldots a_\iota)\iff
\exists z[D'_{x_1\ldots x_\iota y}(a_1\ldots a_\iota z)\and B]$

\medskip

Freeze the variables in $(/ 48)$. By the variant theorem, we can then replace~$D'$
by~D, so we need to prove a closed formula of the form

\nok// 49. $B_z(fa_1\ldots a_\iota)\iff
\exists z[D_{x_1\ldots x_\iota y}(a_1\ldots a_\iota z)\and B]$

\noin (now with variable-free terms $a_1$, \dots, $a_\iota$). Then

\nok// 50. $D_{x_1\ldots x_\iota z}(a_1\ldots a_\iota fa_1\ldots a_\iota)\and
B_z(fa_1\ldots a_\iota)\imp\exists z[D_{x_1\ldots x_\iota y}(a_1\ldots a_\iota z)\and B]$

\noin is a substitution formula, but the first conjunct holds by the defining axiom
and an identity formula, so we have the forward direction of~$(/ 49)$. For the
the backward direction, let r~be the special constant for
$\exists z[D_{x_1\ldots x_\iota y}(a_1\ldots a_\iota z)\andD B]$.
By the special axiom for~r and a closed instance of~EC, we then have

\nok// 51. $D_{x_1\ldots x_\iota y}(a_1\ldots a_\iota r) \and B_z(r)$

\noin and we have

\nok// 52. $r=fa_1\ldots a_\iota$

\noin by a closed instance of UC, and thus we have the left hand side of $(/ 49)$.

The last statement holds since $T_1$ is a conservative extension of~T.\bul

\medbreak

This treatment of f-extensions is much simpler than anything in the literature, so far as I know.
Shoenfield makes the conservativity of an f-extension, and the properties of the translation,
depend on Herbrand's theorem, and so invokes a superexponential algorithm.

For a defining axiom of the form

\noj $fx_1\ldots x_\iota = y \iff y=b$

\noin UC and EC are trivial, and we write the defining axiom simply as

\noj $fx_1\ldots x_\iota = b$

\noin and call it an \ii explicit definition..

An \ii extension by definitions.\/ of T is a theory $\rm T_1$ obtained by a finite
sequence of \hbox{p-extensions} and f-extensions.
It is possible to join all extensions by definitions of~T into a single theory schema.
Let $\dot T$~be T with all defining axioms adjoined. Then we can write $\vdash_{\dot T}A$
to avoid the cumbersome phrase ``A is a theorem of an extension by definitions of~T''.

This is one advantage of the $p^{\phantom0}_D$ and
$f^{\phantom0}_{\exists yD}$ notation:
it makes it impossible to introduce the same symbol with two different defining axioms.
Another is that in the arithmetization it will suffice to arithmetize just the
symbols of~T; then those of~$\dot T$ can be arithmetized in one benign swoop simply
by recursively replacing the symbols of~T in the subscripts by their arithmetizations.

\section//   12.{Default formulas}
A \ii default formula.\/ is $\Neg\exists xB\impP r=0$, where r is the special constant
for~$\exists xB$. They are true in the structure~$\hat\sigma$ of~$\S/   7$, and now
we show that they can be used in proofs of plain formulas.

\Th//  15. {\it If\/ \ro A is a plain formula and its closure is a tautological
consequence of formulas in\/~$\Delta(T)$ and of default formulas, then\/
$\vdash_T A$.}

\pf Let $u\mapsto u^+$ be defined by replacing, everywhere it appears in~u,
each special constant r for $\exists xB$ by the special constant~$\bar r$ for

\nok// 53. $\exists x\{B \or [\Neg \exists x' B_x(x') \and x=0]\}$

\noin Then $(/ 53)$ is a logical theorem, so

\goodbreak

\nok// 54. $ \exists xB \impP B(\bar r)$

\nok// 55. $ \Neg\exists xB \impP \bar r=0$

\noin are logical theorems. The map $C\mapsto C^+$ is a homomorphism. If C is
is a special axiom, $C^+$ is of the form~$(/ 54)$, and if C~is a default formula,
it is of the form~$(/ 55)$. If C~is a substitution, identity, or equality formula,
so is~$C^+$. Finally, $\bar A^+$ is~$\bar A$ since A~is plain. \bul

\section//   13.{Relativization}
Let $ A_\phi$ be the formula
obtained from A by replacing each $\exists xB$ occurring in it
by $\exists x[\phi(x)\andD B]$.

We use $\phi(\free A)$ as an abbreviation for

\noj $\phi(x_1)\andD\cdots\andD\phi(x_\iota)$

\noin where $x_1$, \dots, $ x_\iota$ are the free variables of A. Write

\nok// 56. $A^\phi$ \quad for\quad $\phi(\free A)\impP A_\phi$

\noin (We adopt the convention that if A is closed, so that $\free A$ is the empty
expression, then $\phi(\free A)\impP A_\phi$ is $A_\phi$.)
We call $ A^\phi$ the \ii relativization of~\ro A by.\/~$\phi$.

Let T~be a theory. We say $\phi$ \ii respects~\ro f in.\/ T in case

\noj $\vdash_T\ \phi(x_1)\andD\cdots\andD\phi(x_\iota)
\imp \phi(\ro fx_1\ldots x_\iota)$

\noin where $\rm I(f)=\iota$, and $\phi$ \ii respects~\ro b in.\/~T in case

\noj $\vdash_T\ \phi(x_1) \andD\cdots\andD\phi(x_\iota)\imp\phi(b)$

\noin where $ x_1$, \dots, $ x_\iota$ are the variables occurring in~b, and that
$\phi$ \ii respects \ro A in.\/ T in case

\noj $\vdash_T A^\phi$

\Th//  16. {\it If for each\/ \ro f occurring in\/ \ro b, $\phi$ respects\/ \ro
f in\/~\ro T, then $\phi$ respects\/~\ro b in\/~\ro T.}

\pf By induction on the formation of terms. \bul

\medskip

We say that $\phi$ is a \ii relativizer of.\/ T in case $\phi$ respects in T each
function symbol of L(T) and each nonlogical axiom of~T. If T~is open, the latter
condition holds automatically.

\Th//  17. {\it Let\/ $\phi$ be a relativizer of\/ \ro T.
Let\/~$A$ be a plain formula of\/ $L(T)$ such that\/
$\vdash_T A^\phi$, and let\/~$T'$ be $T[A]$. If\/
$\vdash_{T'}B$ then\/ $\vdash_T B^\phi$.}

\pf Let $\pi\vdash_{T'}B$.
Replace each C in $\pi$ by $C^\phi$. If C~is a nonlogical
axiom of~$T'$, then $C^\phi$ is a theorem of~T by hypothesis. The mapping
$C\mapsto C^\phi$ preserves tautological consequence. If C~is a special axiom,
identity formula, or an equality formula, then $C^\phi$~is a tautological consequence
of~C. Hence we need only consider substitution formulas.

We use default formulas, which can be eliminated by Theorem/  15.
Also, by \hbox{Theorem/  2} we assume that $\pi$~is an expression of~L(T).
Let C be a substitution formula, of the form $D_x(b) \impP \exists xD$.
Then $C^\phi$ is

\nok// 57. $D_{\phi x}(b) \imp \exists x[\phi(x) \andD D_\phi]$.

\noin Every constant occurring in the variable-free term b is either
a constant~$e$ of~L(T), in which case $\vdash_T\phi(e)$ by hypothesis, or
it is the special constant~r
for a closed instantiation $\exists yD'_\phi$, in which case we have $\phi(r)$
by the special axiom if $\exists yD'_\phi$ holds and by the default formula otherwise.
Hence $\vdash_T\phi(b)$ by \hbox{Theorem/  16}, and consequently $(/ 57)$ is a
theorem of~T.~\bul

\medskip

Such a theory $T[A]$, where for some relativizer $\phi$ of T
the formula $A^\phi$ is a plain theorem of~T, is an \ii r-extension.\/ of~T.
An r-extension $T'$ of~T is not in general a conservative extension
(if it is, then there is no point in performing the relativization). But if
T~is consistent, then so is $T'$ by Theorem/  17, since $[0\ne0]^\phi$ is~$0\ne 0$.

A relativization is the particular case of an interpretation [Sh~\S4.7]
in which a theory is interpreted in itself with the interpretation of each
nonlogical symbol being the symbol itself.
In our proof of the inconsistency of~\PA,
relativizations are the only kind of interpretation that will be used.

An \ii extension by definitions and relativizations., or \ii R-extension.\/ for short,
of~T is a theory $T_\nu$ for which there is a sequence $T_0$, \dots,~$T_\nu$ where
$T_0$ is~T and each $T_\mu$ for $1\le\mu\le\nu$ is a t-, p-, f-, or r-extension of
$T_{\mu-1}$. If A is a formula of $T_\mu$, its \ii image.\/ $* A$ is defined
recursively as follows. If $T_\mu$ is a t-extension, $* A$ is A; if it is a p-extension
or f-extension, $* A$ is $* A^*$ where $A^*$ is the translation of A into $T_{\mu-1}$;
if it is an r-extension, $* A$ is $* A^\phi$ where $\phi$ is the relativizer
of $T_{\mu-1}$. If $\vdash_{T_\nu}A$, then $\vdash_T* A$. If T is consistent, so is
$T_\nu$ since $*[0\ne0]$ is $0\ne0$. If $\pi$~is a proof over~T of~A, then
the \ii reduction.\/ $\star\pi$ is the corresponding proof in~T of~$* A$.
A proof in an R-extension of~T is called a proof \ii over.\/~T.

An R-extension~$\rm T'$ of T is said to be \ii relativizable in.\/ T. A theory
schema~$\rm T^*$ is \ii locally relativizable in.\/~T in case each theory in~$\rm T^*$
is relativizable in~T.

\bigskip

\hfill\vbox{\baselineskip=10pt\halign{\ninerm#\hfil\cr
Blind unbelief is sure to err\cr
And scan His works in vain;\cr
God is His own interpreter\cr
And He will make it plain.\cr}}\smallskip

\hfill\hbox{--- \ninerm W{\sevenrm ILLIAM} \ninerm C{\sevenrm OWPER}}

\goodbreak

\section//   14.{Bounded formulas}
Let L be a language containing a binary predicate symbol which we denote
by $\le$, though it can be arbitrary.
If x does not occur in b, we use the abbreviation

\noj $\exists x \leE b\,A \for \exists x[x\le b \andD A]$

\noin (If we unabbreviate $\exists x[x\le b \andD A]$ we obtain
$\exists x\Neg{\vee}\Neg{\le}xb\Neg A$, so $\exists x\leE b$ abbreviates
$\exists x\Neg{\vee}\Neg{\le}xb\Neg$.)
We call $\exists x \leE b$ a \ii bounded quantifier.\/ (for $\le$).
A formula is \ii bounded.\/ (for~$\le$) in case every occurrence in it of~$\exists x$
is subformula of some $\exists x\leE b$ where x does not occur in~b.
If x~does not occur in~b, we also use the
abbreviation $\forall x \leE b\,A$ for
$\Neg\exists x \leE b\Neg\,A$. This is
bounded and is logically equivalent to $\forall x[x\le b\impP A]$.

\Th//  18. {\it Let\/ $\phi$ respect\/ \ro T and let\/~\ro A be a bounded formula
of\/ $L(\dot T)$. Then}

\nok// 58. $\vdash_T\ \phi(\free A) \imp [A \iffF A_\phi]$

\noin {\it Consequently,}

\nok// 59. $\vdash_T\ [\phi(\free A) \impP A] \iff A^\phi$

\pf The proof of $(/ 58)$ is by induction on the construction of bounded formulas.
If A is atomic, $ A_{\phi}$ is~A. Since $\phi(\free \Neg B)$
is $\phi(\free B)$ and $ [\Neg B]_{\phi}$ is $\Neg[B_{\phi}]$, the result
holds by the induction hypothesis if A is~$\Neg B$. Similarly, since

\noj $\phi(\free [B\orR C])\iff\phi(\free B) \and \phi(\free C)$

\noin is a tautology
and since $[B\orR C]_{\phi}$ is \hbox{$ B_{\phi} \orR C_{\phi}$},
the result holds by the induction hypothesis if A is~$ B\orR C$.
Now suppose that A is the formula
$\exists x[x\le b \andD B]$ where x does not occur in~b.
Then every variable in~b occurs free in~A, so by hypothesis and Theorem/  16,
$\phi$~respects~b. Hence

\noj $\vdash_T \phi(\free A) \imp \{ \exists x[x\le b \andD B]\iff
\exists x[\phi(x)\andD x\le b \andD B]\}$

\noin from which the result follows by the induction hypothesis. This proves $(/ 58)$,
and $(/ 59)$ is a tautological consequence of it by the definition $(/ 56)$ of $A^\phi$.
\bul

\medskip

Call a p-extension of T \ii bounded.\/ in case the definiens is bounded, and an
f-extension \ii bounded.\/ in case EC is bounded.
A \ii bounded extension by definitions.\/
of~T is a theory obtained from~T by a finite sequence of bounded p-extensions and
f-extensions.

\goodbreak

\Th//  19. {\it Let\/ $T'$ be a bounded extension by definitions of\/~$T$, and let\/
$A\mapsto A^*$ be the translation from\/ $L(T')$ to\/~$L(T)$. Then for every bounded
formula\/~$A$ of\/~$L(T')$, $A^*$ is a bounded formula of\/~$L(T)$.}

\pf The proof is by induction on the number of p-extensions and f-extensions.
For each p-extension, the result is immediate from the definition of~$A^*$
in~$\S/   10$. Consider an f-extension as in~$\S/   11$. Since EC is bounded,
the definiens~D is of the form $y\preceq b\andD C$ where y~does not occur in~b and
C~is bounded. Then the result holds by~$(/ 47)$ of $\S/   11$. \bul

\section//   15.{The consistency theorem}
We say that $A_1$, \dots, $A_\nu$ is a \ii special sequence.\/ in case

\noj $\Neg A_1 \orR \cdots \orR \Neg A_\nu$

\noin is a tautology; that is, in case for every truth valuation V we have
$V(A_{\mu})=\F$ for some~$\mu$ with $1\le\mu\le\nu$.

\Th//  20. {\it Let\/ \ro T be an inconsistent theory. Then there is a
special sequence of formulas in~$\Delta(T)$.}

\pf Since T is inconsistent, $\vdash_T 0\ne 0$.
Hence $0\ne 0$ is a tautological consequence of formulas in~$\Delta(T)$.
Since $0=0$ is in~$\Delta(T)$, the result follows. \bul

\medskip

A formula \ii belongs to.\/~r in case it is the special axiom for~r or is a
substitution formula $A_x(a)\impP\exists xA$ where r~is the special constant
for~$\rm\exists xA$. The \ii un-nested rank.\/ of~A is the number of occurrences of~$\exists$ in it;
the \ii nested rank.\/ of~A is the maximal rank of instantiations occurring in~A.
We use \ii rank.\/ as a synonym for nested rank; this terminology differs from that of~[Sh].
The \ii rank.\/ of a special constant is the rank of its subscript, so it is
at least~1.
A formula is \ii in.\/ $\Delta_\rho(T)$ in case it is in~$\Delta(T)$
and does not belong to some special constant of rank strictly greater than~$\rho$.
Then every closed theorem of~T is a tautological
consequence of formulas in~$\Delta(T)$, and so of formulas in~$\Delta_\rho(T)$
for some~$\rho$.

The \ii nested rank.\/ of a formula is the maximal rank of instantiations occurring in it.`
A \ii $\rho$-theory.\/ is a theory in which the nested rank of the nonlogical axioms is
most~$\rho$. A \ii $\rho$-proof.\/ in~T is a proof in~T such that all of its formulas are in~$\Delta_\rho(T)$.
A \ii $(\rho,\lambda)$-proof in T is a $\rho$-proof in T in which all the special constants are of level
at most~$\lambda$.
A theory is \ii $\rho$-consistent.\/ in case it is a $\rho$-theory and
there is no $\rho$-proof in~T of $0\ne0$; otherwise it is
\ii $\rho$-inconsistent.. A synonym for 0-consistent is \ii openly consistent..

The Hilbert-Ackermann consistency theorem is an algorithm for the elimination of quantifiers. We follow
the account if [Sh~\S4.3]. Shoenfield treats the case of open theories (i.e., 0-theories); we consider
$\rho_0$-theories, but the proof is the same.

A special sequence is a \ii $(\rho,\lambda,\kappa)$-special sequence.\/ in case
the maximal rank of special constants occurring in its formulas is at most~$\rho$,
their maximal level is at most~$\lambda$,
and the maximal level of any special constant for a
formula of rank~$\rho$ in the special sequence is of level at most~$\kappa$
(so $\kappa\le\lambda$). It is a \ii $(\rho,\lambda)$--special sequence.\/ in case
it is a $(\rho,\lambda,\kappa)$-special sequence for some~$\kappa$.
Note that a $(\rho,\lambda,0)$-special sequence is a $(\rho-1,\lambda,\lambda)$-special sequence.

\smallskip

\Th//  21. {\it Let\/ \ro T be a $\rho_0$-theory and suppose that there is a
$(\rho,\lambda,\kappa)$-special sequence
of $\nu$~formulas in $\Delta_\rho(T)$ with $\rho>\rho_0$ and $\kappa>0$.
Then there exists a $(\rho,2\lambda,\kappa-1)$-special sequence of at
most $\nu^2$ formulas in~$\Delta_\rho(T)$.}

\pf Let $A_1$, \dots,~$A_\nu$ be such a special sequence.
Let R be the set of special constants of rank $\rho$ for which there is a formula
in the special sequence belonging to it; let M~be the elements of~R of
level~$\kappa$, and let
$S=R\setminus M$. Then our special sequence consists of formulas of the form

\nok// 60. $A,\ \exists xB\impP B_x(r),\ B_x(a)\impP\exists xB$

\noin where each A is in $\Delta_{\rho-1}(T)$ or belongs to an element of~S, and
each~$\exists xB$ is the subscript of an element of~M.
Note that there may be many variable-free terms~a such that \hbox{$B_x(a)\impP\exists xB$}
belongs to~r and is in~$(/ 60)$, where r~in~M is the special constant for $\exists xB$.

I claim that $\exists xB$, where its special constant r is in M,
does not occur in any formula C of $(/ 60)$ other than those belonging to~r. This is immediate
if C~is an identity or equality formula for then it is open, and it is immediate for a closed instance of a
nonlogical axiom, for then every instantiation in it is of rank at most~$\rho_0$,
so suppose that C~is either the special axiom $\exists yD\impP D_y(r')$ or is $D_y(b)\impP\exists yD$.
Since the rank of~r is~$\rho$ and the rank of~$r'$
is at most~$\rho$, $\exists xB$ does not occur in $D_y(r')$ or $D_y(b)$. For the same
reason, it cannot occur in $\exists yD$ unless they are the same formula, which is
impossible since $r'$~is not~r. This proves the claim.

In the special sequence $(/ 60)$, replace each occurrence of each
$\exists xB$ by $B_x(r)$.
Since this mapping preserves tautological consequence, the new sequence is a special
sequence. The formulas~A are unaffected, by the claim.
Each special axiom $\exists xB\impP B_x(r)$ becomes the tautology
$B_x(r)\impP B_x(r)$, so it may be deleted. We obtain the special sequence consisting
of formulas of the form

\nok// 61. $A,\ B_x(a)\impP B_x(r)$

\medskip

Call $(r,a)$ a \ii special pair.\/ for the given special sequence in case r is in~M
and the formula \hbox{$B_x(a)\impP B_x(r)$} is in~$(/ 61)$.
Let $u^{(r,a)}$ be obtained from u by replacing r everywhere it appears in~u
by~a. Note that $B_x(r)^{(r,a)}$
is $B_x(a)$ since r~does not appear in~B. Applying this mapping, which
preserves tautological consequence, to the
special sequence~$(/ 61)$, for each special pair~$(r,a)$ we obtain the special sequence
of formulas of the form

\nok// 62. $A^{(r,a)},\ B_x(a')^{(r,a)}\impP B_x(a)$

\medskip

I claim that the

\nok// 63. $A,\ A^{(r,a)}$

\noin are a special sequence.
For if not, there is a truth valuation~V assigning~\T\ to all these formulas.
Since $(/ 62)$ is a special sequence, for each~$(r,a)$ there exists~$a'$ such that

\noj $V\big(B_x(a')^{(r,a)}\impP B_x(a)\big)=\F$

\noin Therefore for all $(r,a)$ we have $V\big(B_x(a)\big)=\F$.
Hence for all~$(r,a)$, \hbox{$V\big(B_x(a)\impP B_x(r)\big)=\T$},
so V assigns~\T\ to every formula in~$(/ 61)$, which is impossible. This proves the claim.
Note that there are at most~$\nu^2$ formulas in~$(/ 63)$.

We must show that each $A^{(r,a)}$ either is in $\Delta_{\rho-1}(T)$
or belongs to some~$r'$ in~S. If A~is an identity or equality formula or is a
closed instance of a nonlogical axiom, then it is in~$\Delta_{\rho_0}$, and so is~$A^{(r,a)}$.

Now suppose that
A~belongs to~$r'$, where $r'\in S$ or the rank of~$r'$ is at most $\rho-1$.
Then $A^{(r,a)}$ is

\noj $\exists yC^{(r,a)}\impP C^{(r,a)}_y\big(r'^{(r,a)}\big)$

\noin or

\noj $C^{(r,a)}_y\big(b^{(r,a)}\big)\impP\exists yC^{(r,a)}$

\noin Now $r'$ is not r, so $r'^{(r,a)}$ is the special constant for
$\exists yC^{(r,a)}$.
Hence $A^{(r,a)}$ belongs to $r'^{(r,a)}$. Since $r'$
and $r'^{(r,a)}$ have the
same rank, if A~is in $\Delta_{\rho-1}(T)$ then so is $A^{(r,a)}$.
Finally, if $r'\in S$ then r does not appear in $\exists yC$, and hence not in~$r'$, since
its level is at least the level of~$r'$. Consequently,
$r'^{(r,a)}$ is~$r'$, and $A^{(r,a)}$
belongs to the element~$r'$ of~S.

Consider a sequence of special constants $r_1$, \dots, $r_\alpha$ where $r_1$ is r and each $r_\beta$
appears in $r_{\beta+1}$ for all~$\beta$ with $1\le\beta<\alpha$. Then $\alpha\le\lambda$. Now
r~is replaced by~a everywhere it appears in~$A^{(r,a)}$ and each special constant appearing in~a is of
level at most~$\lambda$. Hence the special constants appearing in the formulas of~$(/ 63)$ are of level
at most~$2\lambda$, so it is a $(\rho,2\lambda,\kappa-1)$-special sequence. \bul

\medskip

Call $\alpha$ an \ii $\epsilon_0$-number.\/ in case for all $\beta$, the recursion for~$\beta^\alpha$
terminates. Then, recursively, call~$\alpha$ an \ii $\epsilon_\gamma$-number.\/ in case for all
$\epsilon_{\gamma-1}$-numbers $\delta$,
the recursion for~$\delta^\alpha$ terminates and is an $\epsilon_0$-number.

\medskip

\Th//  22. {\it Let\/ \ro T be a $\rho_0$-theory and suppose that there is a
$(\rho,\lambda)$-special sequence
of $\nu$~formulas in $\Delta_\rho(T)$ with $\rho>\rho_0$, and that $\lambda$ is an $\epsilon_1$-number.
Then there exists a $(\rho-1,2^{2^\lambda})$-special sequence of at
most $\nu^{2^\lambda}$ formulas in~$\Delta_{\rho-1}(T)$.}

\medskip

\pf Iterating $\lambda$ times the result of Theorem/  21,
we obtain a $(\rho-1,2^\lambda \lambda)$-special sequence of at
most $\nu^{2^\lambda}$
formulas. Since \hbox{$2^{2^\lambda}>2^\lambda \lambda$}, it is a $(\rho-1,2^{2^\lambda})$-special sequence. \bul

\medskip

Write $\alpha^\beta$ as $\alpha\uparrow\beta$ and $2\uparrow\cdots\uparrow2\uparrow\lambda$ with $\mu$
occurrence of~2 as $\lambda_\mu$.

\medskip

\Th//  23. {\it Let\/ \ro T be a $\rho_0$-theory and suppose that there is a
$(\rho,\lambda)$-special sequence
of $\nu$~formulas in $\Delta_\rho(T)$ with $\rho>\rho_0$,
and that $\lambda$ is an $\epsilon_{\rho-\rho_0+4}$-number.
Then there exists a special sequence of at
most $(\nu\uparrow2\uparrow\lambda)\uparrow\lambda_{\rho-\rho_0+2}$ formulas in~$\Delta_{\rho_0}(T)$.}

\medbreak

\pf Iterating $\rho-\rho_0$ times the result of Theorem/  22, we find that there is a special sequence
of at most
$$\eqalignno{
&(\cdots((\nu\uparrow2\uparrow\lambda)\uparrow2\uparrow(2\uparrow2\uparrow\lambda))\cdots\uparrow
(2\uparrow\lambda_{\rho-\rho_0}))\cr
=\;&(\nu\uparrow2\uparrow\lambda)\uparrow2\uparrow(2\uparrow2\uparrow\lambda\cdot\ldots\cdot
\lambda_{\rho-\rho_0})\cr
\noalign{\hbox{(using $(\mit x\uparrow \mit y)\uparrow \mit z=\mit x\uparrow(\mit y\cdot \mit z)$)}}
=\;&(\nu\uparrow2\uparrow\lambda)\uparrow2\uparrow2\uparrow2\uparrow(2\uparrow\lambda+\cdots+
\lambda_{\rho-\rho_0-1})\cr
\noalign{\hbox{(using $(2\uparrow \mit y)\cdot(2\uparrow \mit z)=2\uparrow(\mit y+\mit z)$)}}
\le\;&(\nu\uparrow2\uparrow\lambda)\uparrow\lambda_{\rho-\rho_0+2}\cr
}$$
formulas in $\Delta_{\rho_0}(T)$. \bul

\medskip

Under the finitary hypothesis that all numbers are $\epsilon_\gamma$-number for all $\gamma$, this theorem
for $\rho_0=0$ takes the traditional form of the Hilbert-Ackermann consistency theorem: if T is
openly consistent, then it is consistent.

\medskip

\Th//  24. {\it Let\/ \ro T be a $\rho_0$-theory and\/ \ro A a formula of\/ \ro T of nested
rank at most $\rho_0$. If there is a $(\rho,\lambda)$-proof in\/ \ro T of\/ \ro A
and $\lambda$ is an $\epsilon_{\rho-\rho_0+4}$-number, then there is a $\rho_0$-proof
in\/ \ro T of\/ \ro A.}

\pf Apply Theorem/  23 to the inconsistent $\rho_0$-theory $T[\Neg A]$. \bul

\medbreak
\goodbreak

The crucial feature of these theorems is that the bounds involve only rank and level, and not the length
of the proof.

\bigskip

\hfill\vbox{\baselineskip=10pt\halign{\ninerm#\hfil\cr
Do I contradict myself?\cr
Very well, then, I contradict myself,\cr
(I am large---I contain multitudes.)\cr}}\smallskip

\hfill\hbox{--- \ninerm W{\sevenrm ALT} \ninerm W{\sevenrm HITMAN}}

\section//   16.{Simple proofs}
The notion of proof is just what is needed
for theoretical purposes, in particular
for the consistency theorem. But how do we recognize that a
proof is in fact a proof? People often refer dismissively to something as being
a mere tautology, but this is akin to calling the Milky Way a mere galaxy.
It is widely conjectured that $\cal NP\!\ne\!{\rm co\hbox{-}}NP$, implying that proofs of
tautologies may in general be infeasibly long. We need another notion of proof,
such that a mathematician looking at one can see by inspection that it is a proof,
and we need assurance that the task of arithmetizing proofs, and proving that they
are proofs, is a purely mechanical feasible procedure.

Let A be in normal form. A \ii special case.\/ of A is a formula obtained by
applying the following steps zero or more times. If the leftmost occurrence of a
quantifier is in~$\forall xB$, replace $\forall xB$ by~$B_x(b)$ where b~is some variable-free term;
otherwise, if it is in $\exists xB$, replace $\exists xB$ by $B_x(r)$ where r~is the special
constant for $\exists xB$.  If $A'$ is a special case of~A, then $A\impP A'$ is a logical theorem.

An \ii equality substitution.\/ is

\nok// 64. $\rm a=b\andD A_x(a) \impP A_x(b)$

\noin It is a logical theorem by Theorem~$/  4$.

A \ii literal.\/ is a variable-free atomic formula or the negation of a variable-free
atomic formula;
the \ii opposite.\/ of a literal~A is $\rm\Neg A$ if A~is atomic or is~B if~$\rm A$
is~$\rm\Neg B$.

The \ii conjuncts.\/ of A are defined recursively as follows: if A is not a conjunction,
the conjuncts of~A consists of A~alone; the conjuncts of $B\andD C$ are the conjuncts
of~B together with the conjuncts of~C. The \ii disjuncts.\/ of~A are defined similarly.
Call C a \ii proper disjunct.\/ of~A in case it is a disjunct of~A but is not~A.

Say that $D'$ \ii follows by one-resolution.\/ from C and D in case C~is a literal,
its opposite is a proper disjunct of~D, and $D'$ is the result of deleting this disjunct
from~D.

Suppose that we want to prove the implication $\rm A\impP B$. In a direct proof,
we assume the hypothesis~A, perform a chain of reasoning, and deduce~B. In an indirect
proof, we assume A and~$\rm\Neg B$, perform a chain of reasoning, and derive a
contradiction. Indirect proofs are twice as efficient, since we have two starting
points rather than one, and we shall always use indirect proofs.

Now we describe the notion of a \ii very simple proof.\/ of A in the theory T.
It consists of a list of formulas $A_0$, \dots,~$A_\nu$ with the following properties.
$A_0$ is the normal form of $\Neg\bar A$. Each $A_\mu$ for $1\le\mu\le\nu$ is
a conjunct of a special case of a nonlogical axiom, is an equality substitution,
or follows by one-resolution from $A_\alpha$ and~$A_\beta$ for some $\alpha$
and~$\beta$ strictly less than~$\mu$.
Finally, either a literal and its opposite occur in the list or a formula $a\ne a$
occurs in the list. Then

\noj $A_0\impP A_1,\quad \ldots,\quad A_0\impP A_\nu,\quad \bar A$

\noin preceded by the proofs of the equality substitutions and special cases,
is a proof of A.

A \ii simple proof.\/ of A in T consists of a finite sequence $C_1$, \dots, $C_\kappa$
of variable-free formulas, called \ii claims., and very simple proofs of

\noj $C_1,\quad C_1\impP C_2,\quad \ldots,\quad
C_1\andD\cdots\andD C_{\kappa-1}\impP C_\kappa,\quad
C_1\andD\cdots\andD C_\kappa\impP A$
o

\noin Here is a convenient way to organize a simple proof. String all the very
simple proofs together, introduce each claim by assuming its negation, and when
the claim is established delete (or never use again) the negation of the claim
and all the subsequent steps before the establishment of the claim.

\bigskip

\hfill\vbox{\baselineskip=10pt\halign{\ninerm#\hfil\cr
You are at liberty to make any possible Supposition:\cr
And you may destroy one Supposition by another:\cr
But then you may not retain the Consequences,\cr
or any part of the Consequences of the Supposition so destroyed.\cr}}\smallskip

\hfill\hbox{--- \ninerm G{\sevenrm EORGE} \ninerm B{\sevenrm ERKELEY}}

\section//   17.{Register machines}
There are many different versions of Turing machines in the literature.
The register machines of Leivant [Le~\S3.1]\foot{ Daniel Leivant,
``Ramified recurrence and computational complexity I:
Word recurrence and poly-time'', 320-343 in P. Cole and J. Remmel,
eds., {\it Feasible Mathematics II, Perspectives in Computer Science},
Birkhauser-Boston, 1994.}
are particularly appealing: there are no cursors and no tape symbols other than bits,
and the commands are direct expressions of four basic string functions which we here denote by
0~1~P~C: 0 prepends a zero bit; 1 prepends a one bit; P~is the string predecessor, which deletes
the leftmost bit (if any); and C is the case function
$$C(X,Y,Z,W)=\cases{Y&if X is empty\cr
Z&if X begins with 0\cr
W&if X begins with 1\cr}$$

Let $0_\lambda$ be the string $0\ldots0$ with $\lambda$ occurrences of 0, and similarly for $1_\lambda$.
If f~is a unary function symbol, let $f^\lambda$ be $f\ldots f$ with $\lambda$~occurrences of~f.

We use a special case of Leivant's register machines.
By a \ii register machine of index.\/ $\iota$ we mean $\kappa$ \ii states., which we take to be the
strings~$0_\lambda$ for $1\le\lambda\le\kappa$ with 0~the \ii initial state.\/ and $0_\kappa$ the
\ii halt state.,
and \ii registers.\/ $\rho_\nu$ for $1\le\nu\le\mu$, where $\mu>\iota$. The registers~$\rho_1$,
\dots,~$\rho_\iota$ are the \ii input registers.\/ and $\rho_\mu$~is the \ii output register..
Each state~$0_\lambda$ other than
the halt state has a unique \ii command.. The commands can be any of the following, called respectively
0, 1, P, and C commands:
$$0\lambda\nu\nu'\lambda'\qquad
1\lambda\nu\nu'\lambda'\qquad
P\lambda\nu\nu'\lambda'\qquad
C\lambda\nu\lambda_1\lambda_2\lambda_3$$

A \ii configuration.\/ is a $(\mu{+}1)$-tuple $\vec Y$ where $Y_{\mu+1}$ is a state.
We say that $Y_{\nu}$ is \ii stored in the
register.\/~$\rho_{\nu}$. It is an \ii initial configuration.\/ in case $Y_{\mu+1}$ is the initial
state and all registers other than the input registers are empty.

Consider a configuration in which its state $0_\lambda$
is not the halt state. If its command is a 0, 1, or P command,
store $0Y_\nu$, $1Y_\nu$, or $PY_\nu$ respectively in $\rho_{\nu'}$,
store $0_{\lambda'}$ in $\rho_{\mu+1}$, and
leave the strings stored in the other registers unchanged.
If its command is a C~command, store $0_{\lambda_1}$, $0_{\lambda_2}$, or~$0_{\lambda_3}$ in $\rho_{\mu+1}$
according as $Y_\nu$ is empty, begins with~0, or begins with~1, and leave the strings stored in the other
registers unchanged. The result is the \ii next.\/
configuration, obtained in one \ii step.. If the machine is in the halt state, the next configuration
is defined to be unchanged. If for all $\iota$-tuples
$\vec X$ the machine whose initial configuration has $\vec X$ in the input registers
reaches the halt state after a number of steps
bounded by a polynomial in $|\vec X|$, with Y~in the output register, then $\vec X\mapsto Y$
is in PTF, and conversely.

In [Le~\S3.3], Leivant proves that the function

\nok// 65. $(Y_1,\ldots,Y_\mu,Y_{\mu+1})\mapsto(Y_1',\ldots,Y_\mu',Y_{\mu+1}')$

\noin taking one configuration to
the next can be expressed using just the basic functions, without recursion. For $1\le\lambda<\kappa$ and
$1\le\gamma\le\mu+1$, if $0_\lambda$ has a 0~command, let
$$h^\lambda_\gamma(\vec y)=\cases{0y_\nu&if $\gamma=\nu$\cr
\noalign{\vskip1pt}
0^{\lambda'}\epsilon&if $\gamma=\mu+1$\cr
y_\gamma&otherwise\cr}$$
Similarly, if $0_\lambda$ has a 1~command or a P~command, replace $0y_\nu$ by $1y_\nu$ or $Py_\nu$ respectively.
If $0_\lambda$ has a C~command, let
$$h^\lambda_\gamma(\vec y)=\cases{C(y_\nu,0^{\lambda_1}\epsilon,0^{\lambda_2}\epsilon,0^{\lambda_3}\epsilon)&if
$\gamma=\mu+1$\cr
y_\gamma&otherwise\cr}$$
Let $h^\kappa_\gamma(y_\gamma)=y_\gamma$. The $h^\lambda_\gamma$ express
$(/ 65)$ when the state $Y_{\mu+1}$ is known.

Let $C'(w,y,z)$ abbreviate $C(w,y,z,\epsilon)$ and let

\nok// 66. $g_\gamma(\vec y)=C'(Py_{\mu+1},h^1_\gamma(\vec y),C'(P^2y_{\mu+1},h^2_\gamma(\vec y),\ldots,
C'(P^\kappa_\gamma(\vec y),\epsilon)\ldots))$

\noin Then $\vec g$ is $(/ 65)$.

\section//   18.{Incompleteness without diagonalization}
The \ii Kolmogorov complexity.\/ of a number $\xi$, denoted by $\K(\xi)$, is the length
(number of bits in the program) of the
shortest register machine that halts and outputs $\xi$. The notion was first
introduced by Solomonoff%
\yy\foot{R. J. Solomonoff, ``A preliminary report on a general theory of inductive
inference'', {\it Report V-131}, Zator Co., Cambridge, Massachusetts, Feb. 1960, revised
Nov. 1960.
\pdfklink{world.std.com/~rjs/z138.pdf}{world.std.com/~rjs/z138.pdf}}\zz\
and independently by Kolmogorov\foot{A. N. Kolmogorov, ``Three approaches to the quantitative definition of
information'', {\it Problems of Information Transmission}, 1, 1-7, 1965.} (who then acknowledged
Solomonoff's priority).

The most important theorem in the subject was proved by Chaitin.%
\foot{G. J. Chaitin, ``Computational complexity and G\umlaut odel's incompleteness
theorem'', {\it ACM SIGACT News}, 9, pp.~11-12, 1971.}
Let T be a theory capable of arithmetizing itself and expressing combinatorics,
such as~\PA\ or, as we shall see,~\Q.
Let $\bar\xi$~ be the numeral S\dots S0 with $\xi$~occurrences of~S.
Let K~be the arithmetization of~\K. Choose a number~$\lambda$ and
construct a machine as follows. It systematically searches through strings, ordered first by
length and then lexicographically. If it finds one that is a proof in~T
of $K(\bar\xi)>\bar\lambda$,
it outputs~$\xi$ and halts. Now if $\lambda$ is itself of low complexity,
it can be referred to by a term of~T of length much less than~$\lambda$,
so for a large enough such~$\lambda$ the
machine is itself of length less than~$\lambda$. Fix such a~$\lambda$ and denote it
by~$\ell$. Call the corresponding machine the \ii Chaitin machine for.\/~T.

We have

\nok// 67. if $\K(\xi)\le\ell$, there is a $\pi$ such that
$\pi\vdash_TK(\bar\xi)\le\bar\ell$

\noin simply by verifying that the steps of the register machine in question are followed.
Chaitin's theorem is

\nok// 68. if T is consistent, there do not exist $\xi$ and $\pi$ such that
$\pi\vdash_TK(\bar\xi)>\bar\ell$

\noin for otherwise the Chaitin machine would find a proof in T of some
$K(\bar\eta)>\bar\ell$
(where $\eta$ may or may not be the same as $\xi$) and output~$\eta$.
Then we would have $\K(\eta)\le\ell$ by the definition of Kolmogorov complexity,
giving a contradiction in~T by~$(/ 67)$.
G\umlaut odel's first incompleteness theorem is a consequence.
There is no diagonalization or self-reference in this proof;
Chaitin remarks that it is a version of the Berry paradox.

Given the proof of $K(\bar\xi)>\ell$, there are exponentially many strings for the Chaitin machine
to search before finding the proof of $K(\bar\eta)>\ell$, so the proof of Chaitin's theorem just given
assumes that T~not only
proves the consistency of its own arithmetization but proves that exponentiation is total. This
observation is due to Robert Solovay (personal communication).

Kritchman and Raz [KrRa]%
\yy\foot{Shira Kritchman and Ran Raz, ``The surprise examination and the second
incompleteness theorem'', {\it Notices of the AMS}, 57, 1454-1458, 2010.
\pdfklink{www.ams.org/notices/201011/rtx101101454p.pdf}{www.ams.org/notices/201011/rtx101101454p.pdf}}\zz\
have given a stunning proof without self-reference or diagonalization
of G\umlaut odel's second incompleteness theorem.

There are strictly fewer than $2^{\ell+1}$ register machines with at most $\ell$
bits, so
by the pigeonhole principle there is at least one $\xi$ with
$\xi<2^{\ell+1}$ such that $\K(\xi)>\ell$.
Let $\delta$ be the number of such $\xi$, so

\nok// 69. $\delta>0$

\medskip

There are $2^{\ell+1}$ days left in the course and the teacher announces that there
will be an examination on one of those days, but it will come as a surprise. Think
of~$\delta$
as being the number of days remaining after the surprise examination. The examination is
not given on the last day, by~$(/ 69)$. Now
suppose that $\delta=1$ (the surprise examination occurs on the penultimate
day of classes). Then
there is a unique~$\xi$ with $\xi<2^{\ell+1}$ such that $\K(\xi)>\ell$,
and T~proves $K(\bar\eta)\le\bar\ell$ for every other~$\eta$ with
$\eta<2^{\ell+1}$, by~$(/ 67)$.
Hence if $\delta=1$, T~proves $K(\bar\xi)>\bar\ell$, which is impossible if T~is
consistent, by Chaitin's theorem.
Consequently, if~T proves the consistency of its arithmetization then
T~proves $\bar\delta\ge2$. Now suppose that $\delta=2$ and argue in the same way, and
continue up to $\delta=2^{\ell+1}$. In this way, if T~proves the consistency of its
arithmetization, it proves a contradiction, since $\delta\le2^{\ell+1}$.
This yields G\umlaut odel's second incompleteness theorem.

This proof is radically different from G\umlaut odel's self-referential proof.
The latter derives a contradiction from the consideration of proofs
of an entirely unrestricted kind, whereas the former derives a contradiction
from the consideration of proofs of quite specific
kinds, of limited complexity. We shall exploit this difference.

} % end \everymath={\rm} \everydisplay={\rm}

\vfill\eject %%%%%%%%%%%%%%%%%%%%%%%%%%%%%%%%%%%%%%%% REMOVE

\section//   19.{The inconsistency theorem}

\#1.~We construct an open theory \SS.
The nonlogical symbols of \SS\ are the constant $\ep$, the unary function symbols $\z\ \o\ \pr$ and the
binary function symbols $\oplus\ \odot$.
Semantically, the individuals are strings (concatenations of zero bits and one bits) and
$\ep$~is the empty string, $\z$ prepends a zero bit,
$\o$ prepends a one bit, $\pr$ deletes the leftmost bit, if any, and send the empty string to itself,
$\oplus$ concatenates two strings, and $\odot$, the \ii zero product., sends two strings to the string
all of whose bits are the zero bit and whose length is the product of the lengths of the two strings.

\medskip

\#2.~We construct a theory schema \SS* such that \SS* is locally relativizable in \SS. That is,
for each A, if $\rm\pi\vdash_{\SSs^*}A$ then $\rm\star\pi\vdash_{\SSs} *A$ (where the notation is that
of~$\S/   13)$.
\SS* expresses combinatorial constructions for which there are polynomial length bounds on the objects
being constructed. \SS* has an axiom schema of bounded string induction allowing one to perform string induction
$\rm A_x(\ep)\andD\forall x[A\impP\rm A_x(\z\rm x)\andD\rm A_x(\o\rm x)]\impP\rm A$
on formulas~A for which each quantifier has a polynomial length bound.

\medskip

\#3.~\SS* arithmetizes \SS\ as $\ari\SS$ and proves its open consistency.

\medskip

\#4.~Each string X has a \ii name.\/ $\rm\star X$, consisting of $\beta_i\ldots\beta_\nu\ep$ where
$\nu$ is the length of~X and each $\beta_\mu$ is $\z$ or $\o$ according as the $\mu$th bit of X is the
zero bit or the one bit. Arithmetize this in \SS* by the bounded function symbol Name.
The metamathematical function $\rm u\mapsto\ari u$ replaces each symbol in u by the corresponding encoded
symbol and takes the concatenation. Sometimes we wish to pass the values of variables, rather than the
variables themselves, to the encoding. Let $\rm\aari A$ be the encoded closed formula obtained by replacing
each encoded free occurrence of $\rm\ari x$ in $\rm\ari A$ by $\rm Name\ x$.

\medskip

\#5.~\SS* expresses the notion of a machine (register machine of index 0 as in $\S/   17$) by a bounded
unary predicate symbol, and Kolmogorov complexity by an unbounded unary function symbol~K (the value of~K
is bounded, but its definition uses the utterly unbounded notion of a halting machine).

\medskip

\#6.~Use $\rm\pi\vdash_{\SSs^*}^{\rho\lambda}A$ to mean that $\pi$ is a simple $(\rho,\lambda)$-proof
(i.e., of rank at most~$\rho$ and level at most~$\lambda$) in \SS* of~A, and similarly for \SS\ rather than~\SS*.
Without the initial~$\pi$ it means that there is a~$\pi$ such that it holds. Let
$\'encpf'(p,\bar\rho,\bar\lambda,A)$ be the formula of~\SS* expressing that $p$ is an encoded simple
$(\rho,\lambda)$-proof in~$\ari\SS$ of~$A$.

\medskip

\#7.~Let \FI\ abbreviate ``finitary reasoning''. We shall establish various results of the form that \FI\
shows that \SS* proves~A, but we will not incorporate any additional finitary reasoning into~\SS* itself.
Practically every paragraph of this section should be understood as beginning with ``\FI\ shows that \dots''.
Minor use will be made of \FI\ to treat $\rm*A$ and~$\star\pi$, obtained by a sequence, which may be very long,
of translations and relativizations, as formulas and proofs of~\SS, and to treat numbers, which may
be defined as variable-free terms of~\SS*, as numerals~$\bar\nu$, defined as $\z\ldots\z\ep$ with
$\nu$~occurrences of~$\z$.

\medskip

\#8.~Let $\rho_0$ be the rank of $*\ari{\ro K(x)\le k}$. Let $\ro T_\kappa$ be the set of all strings~X
with $\rm|X|\le\kappa$, and for each subset~S of $\ro T_\kappa$ let $\rm A_{\kappa S}$ be the conjunction
of the formulas $\rm*\ari{\ro K(\star X)\le\bar\kappa}$ for $\rm X\in T_\kappa\setminus S$ and of
$\rm\Neg{*}\ari{K(\star X)\le\bar\kappa}$ for $\rm X\in S$. Let $\rm A_{\kappa\delta}$ be the disjunction of
the $\rm A_{\kappa S}$ with $\rm\#S=\delta$, where \#S is the cardinality of~S. By propositional reasoning
(so by a simple $(\rho_0,0)$-proof), \SS\ proves that a unique
such conjunction holds. (Then $\delta$ is the number of days left after the surprise
examination in the Kritchman-Raz proof.)

\medskip

\#9.~There are $\rho_1$ and $\lambda_1$ such that if $\K(\ro X)\le\kappa$ then
$\vdash_{\SSs^*}^{\rho_1\lambda_1}\rm K(\star X)\le\bar\kappa$, so there are $\rho_2$ and $\lambda_2$ such that
if $\K(\ro X)\le\kappa$ then
$\vdash_{\SSs}^{\rho_2\lambda_2}*[\rm K(\star X)\le\bar\kappa]$.
Hence there are $\rho_3$ and $\lambda_3$ such that
this arithmetizes as $\vdash_{\SSs^*}^{\rho_3,\lambda_3} \ro K(x)\le\bar\kappa \impP
\exists p[\'encpf'(p,\bar\rho_2,\bar\lambda_2,*\aari{\ro K(x)\le\bar\kappa}]$.

\medskip

\#10.~There are $\rho_4$ and $\lambda_4$ such that for all $\rho$ and $\lambda$ there is a $\kappa$
such that \SS* constructs the Chaitin machine
$\ro C^{\rho\lambda}_\kappa$ that searches systematically through all encoded names of $\ari\SS$,
ordered first by length and then lexicographically, for an encoded simple $(\rho,\lambda)$-proof
of some $\ari{\Neg}\oplus*\aari{\ro K(x)\le\bar\kappa}$, and if it finds one halts and returns $\'Name ' x$,
and $\vdash_{\SSs^*}^{\rho_4\lambda_4}\ro C\le\bar\kappa$.

\medskip

\#11.~There are $\rho_5$ and $\lambda_5$ such that there is a $(\rho_5,\lambda_5)$-proof in
\SS* of the main lemma (Theorem/  21) for the Hilbert-Ackermann consistency theory in the form:
if $s$ is a special $(r,l,m)$-sequence for $\ari\SS$ with $r>0$ and $m>0$,
then there is a special $(r,l,m-1)$-sequence $s'$ with $s'\preceq s\odot s$.

\medskip

\#12.~There are $\rho_6$ and $\lambda_6$ such that for all sufficiently large
$\rho$ and $\lambda$ there is a simple
$(\rho_6,\lambda_6)$-proof in~\SS* that if the Chaitin machine~$\rm C^{\rho\lambda}_\kappa$ halts
then there is an encoded simple $(\rho,\lambda)$-proof of $\ari{0\ne0}$.
For if it does halt, we have an encoded simple $(\rho,\lambda)$-proof in~$\ari\SS$ of
$\ari\Neg\oplus*\aari{\ro K(x)\le\bar\kappa}$. But then $\ro K(x)\le\bar\kappa$ by definition of
Kolmogorov complexity and~\#10, and by~\#9 there is a simple $(\rho,\lambda)$-proof in~$\ari\SS$ of
$*\aari{\ro K(x)\le\bar\kappa}$, a contradiction.

\medskip

\#13.~\FI\ shows that there are $\rho_7$ and $\lambda_7$ such that for all sufficiently large $\rho$
and $\lambda$ there is a simple $(\rho_7,\lambda_7)$-proof in~\SS* that the Chaitin machine
$\rm C^{\rho\lambda}_\kappa$ does not halt. Use \#12 and \#11. \FI\ shows that a specific (albeit
enormous; see Theorem/  23) number of iterations of~\#12 produces an open contradiction in~$\ari\SS$,
which is impossible by~\#3. Second, \FI\ shows that there are $\rho_8$ and $\lambda_8$ such that for all
sufficiently large $\rho$ and~$\lambda$ there is
a simple $(\rho_8,\lambda_8)$-proof in~\SS\ of $*[\ro C^{\rho\lambda}_\kappa\' does not halt']$.
This is because the bounds $s\odot s$, $(s\odot s)\odot(s\odot s)$, and so forth iterated a certain number
of times (depending on $\rho$ and~$\lambda$) are terms of~\SS, and hence the reduction of the proof to~\SS\
is of fixed rank and level independent of $\rho$ and~$\lambda$.
Third, arithmetizing this result, \FI~shows that there are $\rho_9$ and $\lambda_9$ and $\rho_{10}$ and
$\lambda_{10}$ such that for all sufficiently large $\rho$ and~$\lambda$,
$\vdash_{\SSs^*}^{\rho_9\lambda_9}\exists p
[\'encpf'(p,\bar\rho_{10},\bar\rho_{11},\ari{*\ro C^{\rho\lambda}_\kappa\rm{\ does\ not\ halt}})]$.

\medskip

\#14.~\SS* expresses string exponentiation by an unbounded binary
function symbol~$\uparrow$ with default value $\ep$ if
the recursion does not terminate. (The intended semantics is that $x\uparrow y$ denotes the string all of whose
bits are the zero bit, of length $|x|^{|y|}$, if it exists.)
In \SS* define $\epsilon(y)\iffF\forall x[x\uparrow y\ne0]$.
This is an unbounded unary predicate symbol, expressing that for all~$x$ the recursion for $x\uparrow y$
terminates. We cannot prove $\forall x[\epsilon(x)]$ in \SS*, but for any specific string~X, \SS* proves
$\epsilon(\star\ro X)$. Similarly, there are $\rho_{11}$ and $\lambda_{11}$ and $\rho_{12}$ and $\lambda_{12}$
such that
$\vdash_{\SS^*}^{\rho_{11}\lambda_{11}}\exists p
[\'encpf'(p,\bar\rho_{12},\bar\lambda_{12},*\epsilon({\rm Name}\ x))]$.

\medskip

\#15.~Now we are ready to exploit the Kritchman-Raz proof. Fix large $\rho^*$ and $\lambda^*$ and the Chaitin
machine $\rm C=C^{\rho^*\lambda^*}_\kappa$. Some of the $\rho$'s and $\kappa$'s with subscritps in the
following argument depend on~$\kappa$, and therefore on $\rho^*$ and~$\kappa^*$, but only via iterated
exponentials of~$\kappa$. A machine can be built to construct values of primitive recursive functions, such as
$\rho^*$ and~$\lambda^*$, with the length of the machine growing only linearly with the length of the equations
defining the primitive recursive functions (not their values), so we choose $\rho^*$ and~$\lambda^*$ sufficiently
large so that the $\rho$'s and ~$\lambda$'s with subscripts will be less than $\rho^*$ and~$\lambda^*$
(much less, for a comfort zone). We shall stop making these $\rho$'s and $\lambda$'s explicit.

\medskip

\#16.~Recall the formulas $\rm A_{\kappa\delta}$ of \SS\ from \#8. \SS* proves $\rm\Neg A_{\kappa0}$ by
the pigeonhole principle. Now argue as follows. Suppose that $\rm A_{\kappa1}$. Then there is a unique~X
in~$\rm T_\kappa$ such that for all other~Y in~$\rm T_\kappa$ we have $\rm\K(Y)\le\kappa$. By~\#9,
there are proofs in~\SS* of $\rm K(\star Y)\le\bar\kappa$ for all of them and, since $\rm\Neg A_{\kappa0}$,
a proof in~\SS* of $\rm\Neg K(\star X)$. Consequently, there is a~$\pi$ such that
$\rm\pi\vdash_{\SSs}\Neg{*}K(\star X\le\bar\kappa)$. If we knew that~$\pi$ were exponentiable, we
would have a bound on the number of encoded names the Chaitin machine would have to search to find
a simple encoded proof of some $\rm\ari{\Neg{*}K(\star X')\le\bar\kappa}$, which is impossible
by~\#13. But by~\#14, we have a proof in~\SS* of
$\exists p[\'encpf'(p,\bar\rho_{12},\bar\lambda_{12},*\epsilon(\'Name '\ari\pi)]$ and so an encoded proof
that C~halts, which is impossible by the third part of~\#13. Hence $\rm\Neg A_{\kappa1}$.
Proceed in this way for all the $\rm S_{\kappa\delta}$. Thus \FI\ shows that \SS*, and consequently~\SS,
is inconsistent.

\medskip

\#17.~\FI\ can be expressed in \PA, so \PA\ proves that \SS\ (arithmetized in \PA) is inconsistent.
But \PA\ proves the full Hilbert-Ackermann consistency theorem, so \PA\ also proves that \SS\ is
consistent. Hence Peano Arithmetic is inconsistent.

\medskip

Here is a summary of the argument. The starting point is the unusual feature of \SS\ that a theory schema \SS*
locally relativizable in it
proves the open consistency of its arithmetization $\ari\SS$.
We modify the Chaitin machine so that it searches for a simple encoded
$(\rho^*,\lambda^*)$-proof in $\ari\SS$ of $\aari{*\ro K(x)>\bar\kappa]}$ for some $x$, where $\rho^*$ and
$\lambda^*$ are chosen so that finitary reasoning shows that
the Kritchman-Raz proof can be carried out by arguments of lower complexity.

The entire proof that Peano Arithmetic is inconsistent is far less intricate than
many proofs in contemporary mathematics.
The inordinate length of the sequel is in deference to the well-known dictum:

\bigskip

\vbox{
\hfill\vbox{\baselineskip=10pt\halign{\ninerm#\hfil\cr
Extraordinary claims require extraordinary evidence.\cr}}\smallskip

\hfill\hbox{--- \ninerm C{\sevenrm ARL} \ninerm S{\sevenrm AGAN}}

}

\vfill\eject

\chap3.{String Arithmetic}
%\" \dd//a. \z x = 0 \"
%\" \dd//b. \o x = 0 \"
%\" \dd//c. \pr x = 0 \"
%\" \dd//d. \ep = 0 \"
%\" \dd//e. x {\oplus y} = 0 \"
%\" \dd//f. x \odot y = 0 \"

\section//   20.{The theory \SS}

\" \a//1. \phantom0 \z x \ne \ep \"

\" \a//2. \phantom0 \o x \ne \ep \"

\" \a//3. \phantom0 \z x = \z y \imp x = y \"

\" \a//4. \phantom0 \o x = \o y \imp x = y \"

\" \a//5. \phantom0 \z x \ne \o y \"

\" \a//6. \phantom0 \ep \oplus x = x \"

\" \a//7. \phantom0 x \oplus \ep = x \"

\" \a//8. \phantom0 \z ( x \oplus y ) = \z x \oplus y \"

\" \a//9. \phantom0 \o ( x \oplus y ) = \o x \oplus y \"

\" \a//10. x \oplus ( y \oplus z ) = ( x \oplus y ) \oplus z \"

\" \a//11. \ep \odot y = \ep \"

\" \a//12. \z \ep \odot \z \ep = \z \ep \"

\" \a//13. \z x \odot y = ( \z \ep \odot y ) \oplus ( x \odot y ) \"

\" \a//14. \o x \odot y = \z x \odot y \"

\" \a//15. x \odot ( y \odot z ) = ( x \odot y ) \odot z \"

\" \a//16. x \odot ( y \oplus z ) = ( x \odot y ) \oplus ( x \odot z ) \"

\" \a//17. ( \z \ep \odot x ) \oplus ( \z \ep \odot y ) = ( \z \ep \odot y ) \oplus ( \z \ep \odot x ) \"

\" \a//18. x \odot y = y \odot x \"

\" \a//19. x \odot y = \ep \imp x = \ep \or y = \ep \"

\" \a//20. \pr \ep = \ep \"

\" \a//21. x \ne \ep \imp x = \z \pr x \or x = \o \pr x \"

\" \t//22. \z x \oplus y \ne \ep \"

\sam22.
\"\p/22.
/H : x : y \
/8 ; x ; y \
/1 ; x \oplus y \
\"

\" \t//23. \o x \oplus y \ne \ep \"

\sam23.
\"\p/23.
/H : x : y \
/9 ; x ; y \
/2 ; x \oplus y \
\"

\" \t//24. x \oplus y = \ep \imp x = \ep \and y = \ep \"

\sam24.
\"\p/24.
/H : x : y \
/21 ; x \
/22 ; \pr x ; y \
/23 ; \pr x ; y \
/6 ; y \
\"

\" \d//25. \'Z' x = \z \ep \odot x \"

\" \t//26. \'Z' x \oplus \'Z' y = \'Z' y \oplus \'Z' x \"

\sam26.
\"\p/26.
/H : x : y \
/25 ; x \
/25 ; y \
/17 ; x ; y \
\"

\" \d//27. x \preceq y \iff \exists z [ \'Z' z \oplus \'Z' x = \'Z' y ] \"

\" \t//28. y \odot \ep = \ep \"

\sam28.
\"\p/28.
/H : y \
/11 ; y \
/18 ; \ep ; y \
\"

\" \t//29. \'Z' \ep = \ep \"

\sam29.
\"\p/29.
/H \
/25 ; \ep \
/28 ; \z \ep \
\"

\" \t//30. \ep \preceq x \"

\sam30.
\"\p/30.
/H : x \
/29 \
/27\bw ; \ep ; x ; x \
/7 ; \'Z' x \
\"

\" \t//31. x \preceq x \"

\sam31.
\"\p/31.
/H : x \
/25 ; \ep \
/28 ; \z \ep \
/6 ; \'Z' x \
/27\bw ; x ; x ; \ep \
\"

\" \t//32. \'Z' \'Z' x = \'Z' x \"

\sam32.
\"\p/32.
/H ; x \
/25 ; x \
/25 ; \'Z' x \
/15 ; \z \ep ; \z \ep ; x \
/12 \
\"

\" \t//33. x \preceq \'Z' x \"

\sam33.
\"\p/33.
/H : x \
/27\bw ; x ; \'Z' x ; \ep \
/32 ; x \
/29 \
/6 ; \'Z' x \
\"

\" \t//34. \'Z' x \preceq x \"

\sam34.
\"\p/34.
/H : x \
/27\bw ; \'Z' x ; x ; \ep \
/32 ; x \
/29 \
/6 ; \'Z' x \
\"

\" \t//35. \'Z' x = \ep \imp x = \ep \"

\sam35.
\"\p/35.
/H : x \
/25 ; x \
/19 ; \z \ep ; x \
/1 ; \ep \
\"

\" \t//36. x \preceq \ep \imp x = \ep \"

\sam36.
\"\p/36.
/H : x \
/27\fw ; x ; \ep : z \
/29 \
/24 ; \'Z' z ; \'Z' x \
/35 ; x \
\"

\" \t//37. \pr \z x = x \"

\sam37.
\"\p/37.
/H : x \
/1 ; x \
/21 ; \z x \
/5 ; x ; \pr \z x \
/3 ; x ; \pr \z x \
\"

\" \t//38. \pr \o x = x \"

\sam38.
\"\p/38.
/H : x \
/2 ; x \
/21 ; \o x \
/5 ; \pr \o x ; x \
/4 ; x ; \pr \o x \
\"

\" \t//39. \pr ( \z x \oplus y ) = \pr \z x \oplus y \"

\sam39.
\"\p/39.
/H : x : y \
/8 ; x ; y \
/37 ; x \oplus y \
/37 ; x \
\"

\" \t//40. \pr ( \o x \oplus y ) = \pr \o x \oplus y \"

\sam40.
\"\p/40.
/H : x : y \
/9 ; x ; y \
/38 ; x \oplus y \
/38 ; x \
\"

\" \t//41. x \ne \ep \imp \pr ( x \oplus y ) = \pr x \oplus y \"

\sam41.
\"\p/41.
/H : x : y \
/21 ; x \
/39 ; \pr x ; y \
/40 ; \pr x ; y \
/37 ; x \
/38 ; x \
\"

\" \t//42. \pr \'Z' \ep = \'Z' \pr \ep \"

\sam42.
\"\p/42.
/H \
/29 \
/20 \
\"

\" \t//43. \pr \'Z' \z x = \'Z' \pr \z x \"

\sam43.
\"\p/43.
/H : x \
/25 ; \z x \
/18 ; \z \ep ; \z x \
/13 ; x ; \z \ep \
/12 \
/18 ; \z \ep ; x \
/25 ; x \
/39 ; \ep ; \'Z' x \
/37 ; \ep \
/6 ; \'Z' x \
/37 ; x \
\"

\" \t//44. \pr \'Z' \o x = \'Z' \pr \o x \"

\sam44.
\"\p/44.
/H : x \
/25 ; \o x \
/18 ; \z \ep ; \o x \
/14 ; x ; \z \ep \
/13 ; x ; \z \ep \
/12 \
/18 ; \z \ep ; x \
/25 ; x \
/39 ; \ep ; \'Z' x \
/37 ; \ep \
/6 ; \'Z' x \
/38 ; x \
\"

\" \t//45. \pr \'Z' x = \'Z' \pr x \"

\sam45.
\"\p/45.
/H :x \
/42 \
/43 ; \pr x \
/44 ; \pr x \
/21 ; x \
\"

\" \t//46. \z \'Z' x = \'Z' \z x \"

\sam46.
\"\p/46.
/H : x \
/25 ; \z x \
/18 ; \z \ep ; \z x \
/13 ; x ; \z \ep \
/12 \
/18 ; x ; \z \ep \
/25 ; x \
/8 ; \ep; \'Z' x \
/6 ; \'Z' x \
\"

\" \t//47. \'Z' x \oplus \'Z' y = \'Z' ( x \oplus y ) \"

\sam47.
\"\p/47.
/H : x : y \
/25 ; x \
/25 ; y \
/25 ; x \oplus y \
/16 ; \z \ep ; x ; y \
\"

\" \t//48. \'Z' x \odot \'Z' y = \'Z' ( x \odot y ) \"

\sam48.
\"\p/48.
/H : x : y \
/25 ; x \
/25 ; y \
/25 ; x \odot y \
/15 ; \z \ep ; x ; \z \ep \odot y \
/15 ; x ; \z \ep ; y \
/18 ; x ; \z \ep \
/15 ; \z \ep ; x ; y \
/15 ; \z \ep ; \z \ep ; x \odot y \
/12 \
\"

\" \t//49. x \preceq y \imp \pr x \preceq \pr y \"

\sam49.
\"\p/49.
/H : x : y \
/27\fw ; x ; y : z \
/26 ; z ; x \
/41 ; \'Z' x ; \'Z' z \
/45 ; x \
/45 ; y \
/26 ; \pr x ; z \
/27\bw ; \pr x ; \pr y ; z \
/35 ; x \
/20 \
/30 ; \pr y \
\"

\" \t//50. x \preceq y \and y \preceq z \imp x \preceq z \"

\sam50.
\"\p/50.
/H :x : y : z \
/27\fw ; x ; y : u \
/27\fw ; y ; z : v \
/10 ; \'Z' v ; \'Z' u ; \'Z' x \
/47 ; v ; u \
/27\bw ; x ; z ; v \oplus u \
\"

\" \t//51. x \preceq y \imp z \oplus x \preceq z \oplus y \"

\sam51.
\"\p/51.
/H : x : y : z \
/27\fw ; x ; y : u \
/47 ; z ; x \
/47 ; z ; y \
/10 ; \'Z' z ; \'Z' u ; \'Z' x \
/26 ; z ; u \
/10 ; \'Z' u ; \'Z' z ; \'Z' x \
/26 ; z ; x \
/47 ; x ; z \
/27\bw ; z \oplus x ; z \oplus y ; u \
\"

\" \t//52. x \oplus y \preceq y \oplus x \"

\sam52.
\"\p/52.
/H : x : y \
/47 ; x ; y \
/47 ; y ; x \
/26 ; x ; y \
/27\bw ; x \oplus y ; y \oplus x ; \ep \
/29 \
/6 ; \'Z' ( x \oplus y ) \
\"

\" \t//53. x \preceq y \imp x \oplus z \preceq y \oplus z \"

\sam53.
\" \p/53.
/H : x : y : z \
/51 ; x ; y ; z \
/52 ; x ; z \
/52 ; z ; y \
/50 ; x \oplus z ; z \oplus x ; z \oplus y \
/50 ; x \oplus z ; z \oplus y ; y \oplus z \
\"

\" \t//54. x \preceq y \imp z \odot x \preceq z \odot y \"

\sam54.
\"\p/54.
/H : x : y : z \
/27\fw ; x ; y : u \
/48 ; z ; x \
/48 ; z ; y \
/16 ; \'Z' z ; \'Z' u ; \'Z' x \
/48 ; z ; u \
/27\bw ; z \odot x ; z \odot y ; z \odot u \
\"

\" \t//55. x \preceq y \imp x \odot z \preceq y \odot z \"

\sam55.
\"\p/55.
/H : x : y : z \
/54 ; x ; y ; z \
/18 ; z ; x \
/18 ; z ; y \
\"

\" \t//56. x \preceq y \oplus x \"

\sam56.
\"\p/56.
/H : x : y \
/27\bw ; x ; y \oplus x ; y \
/47 ; y ; x \
\"

\" \t//57. x \preceq x \oplus y \"

\sam57.
\"\p/57.
/H : x : y \
/56 ; x ; y \
/52 ; y ; x \
/50 ; x ; y \oplus x ; x \oplus y \
\"

\" \t//58. \pr x \preceq x \"

\sam58.
\"\p/58.
/H : x \
/20 \
/30 ; \ep \
/8 ; \ep ; \pr x \
/6 ; \pr x \
/9 ; \ep ; \pr x \
/21 ; x \
/56 ; \pr x ; \z \ep \
/56 ; \pr x ; \o \ep \
\"

\section//   21.{The relativization schema}

\" \d~//a1. \sind \iff \phi ( \ep ) \and \forall x [ \phi ( x ) \imp \phi ( \z x ) \andD \phi ( \o x ) ] \"

n\" \d~//a2. \phi^0 ( x ) \iff \forall y [ y \preceq x \imp \phi ( y ) ] \"

\" \t~//a3. \sind \imp [ \phi^0 ( x ) \imp \phi ( x ) ] \"

\sam a3.
\"\p/a3.
/H : x \
/a2\fw ; x ; x \
/31 ; x \
\"

\" \t~//a4. \sind \imp \phi^0 ( \ep ) \"

\sam a4.
\"\p/a4.
/H \
/a1\fw \
/a2\bw ; \ep : y \
/36 ; y \
\"

\" \t~//a5. \sind \imp [ \phi^0 ( x ) \imp \phi^0 ( \z x ) ] \"

\sam a5.
\"\p/a5.
/H : x \
/a2\bw ; \z x : y \
/a2\fw ; x ; \pr y \
/49 ; y ; \z x \
/37 ; x \
/21 ; y \
/a1\fw ; \pr y \
\"

\" \t~//a6. \sind \imp [ \phi^0 ( x ) \imp \phi^0 ( \o x ) ] \"

\sam a6.
\"\p/a6.
/H : x \
/a2\bw ; \o x : y \
/a2\fw ; x ; \pr y \
/49 ; y ; \o x \
/21 ; y \
/38 ; x \
/a1\fw ; \pr y \
\"

\" \t~//a7. \sind \imp [ \phi^0 ( x ) \and u \preceq x \imp \phi^0 ( u ) ] \"

\sam a7.
\"\p/a7.
/H : x : u \
/a2\bw ; u : y \
/50 ; y ; u ; x \
/a2\fw ; x ; y \
\"

\" \d~//a8. \phi^1 ( x ) \iff \forall y [ \phi^0 ( y ) \imp \phi^0 ( x \oplus y ) ] \"

\" \t~//a9. \sind \imp [ \phi^1 ( x ) \imp \phi^0 ( x ) ] \"

\sam a9.
\"\p/a9.
/H : x \
/a8\fw ; x ; \ep \
/a4 \
/7 ; x \
\"

\" \t~//a10. \sind \imp \phi^1 ( \ep ) \"

\sam a10.
\"\p/a10.
/H \
/a8\bw ; \ep : y \
/a4 \
/6 ; y \
\"

\" \t~//a11. \sind \imp [ \phi^1 ( x ) \imp \phi^1 ( \z x ) ] \"

\sam a11.
\"\p/a11.
/H : x \
/a8\bw ; \z x : y \
/8 ; x ; y \
/a8\fw ; x ; y \
/a5 ; x \oplus y \
\"

\" \t~//a12. \sind \imp [ \phi^1 ( x ) \imp \phi^1 ( \o x ) ] \"

\sam a12.
\"\p/a12.
/H : x \
/a8\bw ; \o x : y \
/9 ; x ; y \
/a8\fw ; x ; y \
/a6 ; x \oplus y \
\"

\" \t~//a13. \sind \imp [ \phi^1 ( x ) \and u \preceq x \imp \phi^1 ( u ) ] \"

\sam a13.
\"\p/a13.
/H : x : u \
/a8\bw ; u : y \
/50 ; y ; u ; x \
/a8\fw ; x ; y \
/53 ; u ; x ; y \
/a7 ; x \oplus y ; u \oplus y \
\"

\" \t~//a14. \sind \imp [ \phi^1 ( x_1 ) \and \phi^1 ( x_2 ) \imp \phi^1 ( x_1 \oplus x_2 ) ] \"

\sam a14.
\"\p/a14.
/H : x_1 : x_2 \
/a8\bw ; x_1 \oplus x_2 : y \
/10 ; x_1 ; x_2 ; y \
/a8\fw ; x_2 ; y \
/a8\fw ; x_1 ; x_2 \oplus y \
\"

\" \d~//a15. \phi^2 ( x ) \iff \forall y [ \phi^1 ( y ) \imp \phi^1 ( x \odot y ) ] \"

\" \t~//a16. \sind \imp [ \phi^2 ( x ) \imp \phi ( x ) ] \"

\sam a16.
\"\p/a16.
/H : x \
/a15\fw ; x ; \z \ep \
/18 ; x ; \z \ep \
/a10 \
/a11 ; \ep \
/25 ; x \
/33 ; x \
/a13 ; \'Z' x ; x \
/a9 ; x \
/a3 ; x \
\"

\" \t~//a17. \sind \imp \phi^2 ( \ep ) \"

\sam a17.
\"\p/a17.
/H \
/a15\bw ; \ep : y \
/11 ; y \
/a10 \
\"

\" \t~//a18. \sind \imp [ \phi^2 ( x ) \imp \phi^2 ( \z x ) ] \"

\sam a18.
\"\p/a18.
/H : x \
/a15\bw ; \z x : y \
/a15\fw ; x ; y \
/13 ; x ; y \
/25 ; y \
/34 ; y \
/a13 ; y ; \'Z' y \
/a14 ; \'Z' y ; x \odot y \
\"

\" \t~//a19. \sind \imp [ \phi^2 ( x ) \imp \phi^2 ( \o x ) ] \"

\sam a19.
\"\p/a19.
/H : x \
/a15\bw ; \o x : y \
/a15\fw ; x ; y \
/13 ; x ; y \
/14 ; x ; y \
/25 ; y \
/34 ; y \
/a13 ; y ; \'Z' y \
/a14 ; \'Z' y ; x \odot y \
\"

\" \t~//a20. \sind \imp [ \phi^2 ( x ) \and u \preceq x \imp \phi^2 ( u ) ] \"

\sam a20.
\"\p/a20.
/H : x : u \
/a15\bw ; u : y \
/50 ; y ; u ; x \
/a15\fw ; x ; y \
/55 ; u ; x ; y \
/a13 ; x \odot y ; u \odot y \
\"

\" \t~//a21. \sind \imp [ \phi^2 ( x ) \imp \phi^2 ( \pr x ) ] \"

\sam a21.
\"\p/a21.
/H : x \
/58 ; x \
/a20 ; x ; \pr x \
\"

\" \t~//a22. \sind \imp [ \phi^2 ( x_1 ) \and \phi^2 ( x_2 ) \imp \phi^2 ( x_1 \oplus x_2 ) ] \"

\sam a22.
\"\p/a22.
/H : x_1 : x_2 \
/a15\bw ; x_1 \oplus x_2 : y \
/18 ; x_1 \oplus x_2 ; y \
/16 ; y ; x_1 ; x_2 \
/18 ; y ; x_1 \
/18 ; y ; x_2 \
/a15\fw ; x_1 ; y \
/a15\fw ; x_2 ; y \
/a14 ; x_1 \odot y ; x_2 \odot y \
\"

\" \t~//a23. \sind \imp [ \phi^2 ( x_1 ) \and \phi^2 ( x_2 ) \imp \phi^2 ( x_1 \odot x_2 ) ] \"

\sam a23.
\"\p/a23.
/H : x_1 : x_2 \
/a15\bw ; x_1 \odot x_2 : y \
/15 ; x_1 ; x_2 ; y \
/a15\fw ; x_2 ; y \
/a15\fw ; x_1 ; x_2 \odot y \
\"

\" \a~//BSI. \phi ( \ep ) \and \forall x' [ \phi ( x' ) \impP \phi ( \z x' ) \andD \phi ( \o x' ) ]
\imp \phi ( x ) \"

\section//   22.{Combinatorics}

\" \d//59. y \' ends with ' x \iff \exists z \preceqQ y [ y = z \oplus x ] \"

\" \t//60. y = z \oplus x \imp y \' ends with ' x \"

\sam60.
\"\p/60.
/H : y : z : x \
/59\bw ; y ; x ; z \
/57 ; z ; x \
\"

\" \d//61. x \approx y \iff x \preceq y \and y \preceq x \"

\" \t//62. x \approx x \"

\sam62.
\"\p/62.
/H : x \
/61\bw ; x ; x \
/31 ; x \
\"

\" \t//63. x \approx y \imp y \approx x \"

\sam63.
\"\p/63.
/H : x : y \
/61\fw ; x ; y \
/61\bw ; y ; x \
\"

\" \t//64. x \approx y \and y \approx z \imp x \approx z \"

\sam64.
\"\p/64.
/H : x : y : z \
/61\fw ; x ; y \
/61\fw ; y ; z \
/50 ; x ; y ; z \
/50 ; z ; y ; x \
/61\bw ; x ; z \
\"

\" \d~//b1. \phi_{b1} ( x ) \iff \forall y [ x \oplus y = x \imp y = \ep ] \"

\" \t~//b2. \phi_{b1} ( \ep ) \"

\sam b2.
\"\p/b2.
/H \
/b1\bw ; \ep : y \
/24 ; \ep ; y \
\"

\" \t~//b3. \phi_{b1} ( x ) \imp \phi_{b1} ( \z x ) \"

\sam b3.
\"\p/b3.
/H : x \
/b1\bw ; \z x : y \
/b1\fw ; x ; y \
/8 ; x ; y \
/3 ; x \oplus y ; x \
\"

\" \t~//b4. \phi_{b1} ( x ) \imp \phi_{b1} ( \o x ) \"

\sam b4.
\"\p/b4.
/H : x \
/b1\bw ; \o x : y \
/b1\fw ; x ; y \
/9 ; x ; y \
/4 ; x \oplus y ; x \
\"

\" \t~//b5. \phi_{b1} ( x ) \"

\sam b5.
\"\p/b5.
/H : x \
/BSI , /b1 ; x : x' \
/b2 \
/b3 ; x' \
/b4 ; x' \
\"

\" \t//65. x \oplus y = x \imp y = \ep \"

\sam65.
\"\p/65.
/H : x : y \
/b5 ; x \
/b1\fw ; x ; y \
\"

\" \t//66. y \oplus x = x \imp y = \ep \"

\sam66.
\"\p/66.
/H : y : x \
/47 ; y ; x \
/26 ; y ; x \
/65 ; \'Z' x ; \'Z' y \
/35 ; y \
\"

\" \d~c1. \phi_{c1} ( x ) \iff \forall y \forall z [ x \oplus y = z \oplus y \impP x = z ] \"

\" \t~c2. \phi_{c1} ( \ep ) \"

\sam c2.
\"\p/c2.
/H \
/c1\bw ; \ep : y : z \
/6 ; y \
/7 ; y \
/66 ; z ; y \
\"

\" \t~c3. \phi_{c1} ( x ) \imp \phi_{c1} ( \z x ) \"

\sam c3.
\"\p/c3.
/H : x \
/c1\bw ; \z x : y : z \
/c1\fw ; x ; y ; \pr z \
/21 ; z \
/6 ; y \
/66 ; \z x ; y \
/1 ; x \
/9 ; \pr z ; y \
/8 ; x ; y \
/8 ; \pr z ; y \
/5 ; x \oplus y ; \pr z \oplus y \
/3 ; x \oplus y ; \pr z \oplus y \
\"

\" \t~c4. \phi_{c1} ( x ) \imp \phi_{c1} ( \o x ) \"

\sam c4.
\"\p/c4.
/H : x \
/c1\bw ; \o x : y : z \
/c1\fw ; x ; y ; \pr z \
/21 ; z \
/6 ; y \
/66 ; \o x ; y \
/2 ; x \
/8 ; \pr z ; y \
/9 ; x ; y \
/9 ; \pr z ; y \
/5 ; \pr z \oplus y ; x \oplus y \
/4 ; x \oplus y ; \pr z \oplus y \
\"

\" \t~//c5. \phi_{c1} ( x ) \"

\sam c5.
\"\p/c5.
/H : x \
/BSI , /c1 ; x : x' \
/c2 \
/c3 ; x' \
/c4 ; x' \
\"

\" \t//67. x \oplus y = z \oplus y \imp x = z \"

\sam67.
\"\p/67.
/H : x : y : z \
/c5 ; x \
/c1\fw ; x ; y ; z \
\"

\" \d~//d1. \phi_{d1} ( y ) \iff \forall x \forall z [ y \oplus x = y \oplus z \imp x = z ] \"

\" \t~//d2. \phi_{d1} ( \ep ) \"

\sam 2.
\"\p/d2.
/H \
/d1\bw ; \ep : x : z \
/6 ; x \
/6 ; z \
\"

\" \t~//d3. \phi_{d1} ( y ) \imp \phi_{d1} ( \z y ) \"

\sam d3.
\"\p/d3.
/H : y \
/d1\bw ; \z y : x : z \
/8 ; y ; x \
/8 ; y ; z \
/3 ; y \oplus x ; y \oplus z \
/d1\fw ; y ; x ; z \
\"

\" \t~//d4. \phi_{d1} ( y ) \imp \phi_{d1} ( \o y ) \"

\sam d4.
\"\p/d4.
/H : y \
/d1\bw ; \o y : x : z \
/9 ; y ; x \
/9 ; y ; z \
/4 ; y \oplus x ; y \oplus z \
/d1\fw ; y ; x ; z \
\"

\" \t~//d5. \phi_{d1} ( y ) \"

\sam d5.
\"\p/d5.
/H : y \
/BSI , /d1 ; y : x' \
/d2 \
/d3 ; x' \
/d4 ; x' \
\"

\" \t//68. y \oplus x = y \oplus z \imp x = z \"

\sam68.
\"\p/68.
/H : y : x : z \
/d5 ; y \
/d1\fw :y ; x ; z \
\"

\" \d~e1. \phi_{e1} ( y ) \iff \forall x \forall w \forall z [ y \oplus x = w \oplus z \and \'Z' y =
\'Z' w \imp y = w ] \"

\" \t~e2. \phi_{e1} ( \ep ) \"

\sam e2.
\"\p/e2.
/H \
/e1\bw ; \ep : x : w : z \
/29 \
/35 ; w \
\"

\" \t~e3. \phi_{e1} ( y ) \imp \phi_{e1} ( \z y ) \"

\sam e3.
\"\p/e3.
/H : y \
/e1\bw ; \z y : x : w : z \
/e1\fw ; y ; x ; \pr w ; z \
/29 \
/46 ; y \
/1 ; \'Z' y \
/8 ; y ; x \
/5 ; y \oplus x ; \pr w \oplus z \
/21 ; w \
/46 ; \pr w \
/3 ; \'Z' y ; \'Z' \pr w \
/9 ; \pr w ; z \
/8 ; \pr w ; z \
/3 ; y \oplus x ; \pr w \oplus z \
\"

\ifx\wrapperUsed\nil
\def\thisEnd{\end}
\else
\def\thisEnd{\relax}
\fi

\thisEnd

%% file: NelsonarXivElemMacros.tex
% MAGNIFICATION
% \magnification=\magstep1
\mathcode`: = 58
\mathcode`; = 59
% FONTS
\font\bigrm=cmr7 scaled \magstep4
\font\BBB=msbm10 scaled \magstep4
\font\tenmsa=msam10

\font\bbbold=msbm10
\def\bbb#1{\hbox{\bbbold#1}}
\font\ninerm=cmr9
\font\ninebf=cmbx9
\font\nineit=cmti9
\font\ninei=cmmi9
\font\ninett=cmtt9
\font\eighti=cmmi8
\font\eightrm=cmr8
\font\eightbf=cmbx8
\font\eightit=cmti8
\font\eightsl=cmsl8
\font\eighttt=cmtt8
\font\eightsy=cmsy8

\font\ninesy=cmsy9
\font\eightsy=cmsy8
\font\nineit=cmti9
\font\eightit=cmti8
\font\sixrm=cmr6
\font\sixbf=cmbx6
\font\sixi=cmmi6
\font\sixsy=cmsy6

\font\ninesans=cmss9
\font\eightsans=cmss8

\def\ninepoint{\def\rm{\fam0\ninerm}%
   \textfont0=\ninerm \scriptfont0=\sixrm \scriptscriptfont0=\fiverm
   \textfont1=\ninei \scriptfont1=\sixi \scriptscriptfont1=\fivei
   \textfont2=\ninesy \scriptfont2=\sixsy \scriptscriptfont2=\fivesy
   \textfont3=\tenex \scriptfont3=\tenex \scriptscriptfont3=\tenex
   \textfont\itfam=\nineit \def\it{\fam\itfam\nineit}%
   \textfont\bffam=\ninebf \scriptfont\bffam=\sixbf
    \scriptscriptfont\bffam=\fivebf \def\bf{\fam\bffam\ninebf}%
   \normalbaselineskip=11pt
   \setbox\strutbox=\hbox{\vrule height8pt depth 3pt width0pt}%
   \let\big=\ninebig \normalbaselines\rm}
\def\eightpoint{\def\rm{\fam0\eightrm}%
   \textfont0=\eightrm \scriptfont0=\sixrm \scriptscriptfont0=\fiverm
   \textfont1=\eighti \scriptfont1=\sixi \scriptscriptfont1=\fivei
   \textfont2=\eightsy \scriptfont2=\sixsy \scriptscriptfont2=\fivesy
   \textfont3=\tenex \scriptfont3=\tenex \scriptscriptfont3=\tenex
   \textfont\itfam=\eightit \def\it{\fam\itfam\eightit}%
   \textfont\bffam=\eightbf \scriptfont\bffam=\sixbf
    \scriptscriptfont\bffam=\fivebf \def\bf{\fam\bffam\eightbf}%
   \normalbaselineskip=11pt
   \setbox\strutbox=\hbox{\vrule height8pt depth 3pt width0pt}%
   \let\big=\eightbig \normalbaselines\rm}
\catcode`@=11
\def\ninebig#1{{\hbox{$\textfont0=\tenrm\textfont2=\tensy
  \left#1\vbox to7.25pt{}\right.\n@space$}}}
\catcode`@=12

% PDF
% \pdfpagewidth=8.5true in
% \pdfpageheight=11true in
% \pdfhorigin=1true in
% \pdfvorigin=1true in

% from http://insti.physics.sunysb.edu/~siegel/tex.shmtl (modified)

%\def\pdfklink#1#2{%
%	\noindent\pdfstartlink user
%		{/Subtype /Link
%		/Border [ 0 0 0 ]
%		/A << /S /URI /URI (#2) >>}{\rgbo{0 0 1}{#1}}%
%	\pdfendlink
%	}
\def\yy{\catcode`_=12\catcode`/=12\catcode`~=12\catcode`_=12}
\def\zz{\catcode`/=\active\catcode`~=13\catcode`_=8}
\yy
\def\sam#1.{\smallskip\par{\quad}{{\eightit\pdfklink{Proof.}{http://math.princeton.edu/~nelson/proof/#1.pdf}}}\zz\hskip8pt }
\zz

%\yy
%\def\pdfKlink#1#2{%
%	\noindent\pdfstartlink user
%		{/Subtype /Link
%		/Border [ 0 0 0 ]
%		/A << /S /URI /URI (#2) >>}{\rgbo{0 0 1}{\underline{#1}}}%
%	\pdfendlink
%	}
%\zz

\yy\zz

% QEA MACROS
\def\ul{\hbox{\tenmsa\char"70}}
\def\ur{\hbox{\tenmsa\char"71}}
\def\umlaut#1{\accent "7F #1}
\let\doublequote="
\def\"{}
\def\slash{\char'57}
\catcode`/=13
\def/{}
\def\a#1.#2\par{\smallbreak\par\noindent a$#1$.\quad$#2$}
\def\ar#1.#2\par{\smallbreak\par\noindent r$#1$.\quad$#2$}
\def\at#1.#2\par{\smallbreak\par\noindent t$#1$.\quad$#2$}
\def\at#1.#2\par{\smallbreak\par\noindent t$#1$.\quad$#2$}
\def\t#1.#2\par{\smallbreak\par\noindent t$#1$.\quad$#2$}
\def\tu#1.#2\par{\medbreak\par\noindent t$#1$.~({\eightrm UC})\quad$#2$}
\def\te#1.#2\par{\medbreak\par\noindent t$#1$.~({\eightrm EC})\quad$#2$}
\def\d#1.#2\par{\smallbreak\par\noindent d$#1$.\quad$#2$}
\def\de#1.#2\par{\smallbreak\par\noindent e$#1$.\quad$#2$}
\def\p/#1.#2\"{\begingroup\baselineskip=9pt\everymath={\scriptstyle}\def\rm{\eightrm}\def\ {\hskip 10pt plus 2pt minus 10pt}$\!\!#2$\par\endgroup}
\newif\ifproof

\let\acute=\'
\def\'#1'{{}$\rm#1${}}
\def\iff{\allowbreak\hskip8pt plus2pt minus3pt\mathrel{\hbox{$\leftrightarrow$}}\hskip8pt plus2pt minus3pt}
\def\iffF{\allowbreak\leftrightarrow}
\def\imp{\allowbreak\hskip8pt plus2pt minus3pt\mathrel{\hbox{$\rightarrow$}}\hskip8pt plus2pt minus3pt}
\def\impP{\allowbreak\rightarrow}

\def\andD{\allowbreak\mathrel{\hskip1pt\&\hskip1pt}}
\def\and{\allowbreak\hskip7pt plus2pt minus3pt\mathrel{\&}\hskip7pt plus2pt minus3pt}
\def\orR{\allowbreak\kern1pt\lor\kern1.25pt}
\def\or{\allowbreak\hskip7pt plus2pt minus3pt\mathrel{\hbox{$\vee$}}\hskip7pt plus2pt minus3pt}
\def\leE{{\scriptstyle\le}}

\def\preceqQ{{\scriptstyle\preceq}}
\let\dotaccent=\.
\def\.{\cdot}

% REGISTERS
\newcount\firstpage
\newcount\footno
\newcount\secno
\newtoks\chaptername
\newtoks\sectionname
\newwrite\tabcon
\newif\ifindexmode

% FORMATTING
\predisplaypenalty=0
\def\foot#1{\global\advance\footno by 1
  {\baselineskip=9pt
	\frenchspacing
    \setbox\strutbox=\hbox{\vrule height7pt depth2pt width0pt}%
	\def\textindent##1{\indent\llap{##1}}%
	\def\sl{\eightsl}\def\bf{\eightbf}\def\rm{\eightrm}%
	\def\it{\eightit}\def\tt{\eighttt}%
	\eightrm\footnote{$^{\hbox{{\sevenrm\the\footno}}}$}{#1}}}
\def\footnoterule{\kern-3pt \hrule width 2truein height 0pt}
\headline={\Headline}
\footline={\Footline}
\def\Headline{\ifnum\firstpage=\pageno \hfil \else\ifodd\pageno%
   \righthead\else\lefthead\fi\fi}
\def\Footline{\ifnum\firstpage=\pageno \hss\tenrm\folio\hss \else{}\fi}
\newif\ifrighthead
\def\lefthead{\tenrm\folio\hfil{\eightrm CHAPTER \chapno. \hskip2pt\uppercase\expandafter{\the\chaptername}}\hfil}
\def\righthead{\ifrighthead\hfil{\eightrm{\botmark}}%
   \hfil\tenrm\folio\else\lefthead\fi}
\let\chapno=0
\def\chap#1.#2{\chaptername={#2} \firstpage=\pageno
    \def\chapno{#1}
   \hrule width0pt depth1.5cm
   \centerline{\tenrm CHAPTER\ #1}\vskip .45cm
   \centerline{\bigrm #2}\vskip 1cm \ifindexmode\write\tabcon%
   {#1\  #2 \the\count0}\fi}
\def\section//   #1.#2{\goodbreak\rightheadtrue\def\sectionname{\eightrm\uppercase{#2}}\secno=#1\bigbreak\noindent{\bf#1.\  #2}\mark{{\eightrm\S}\the\secno.\ \ \sectionname}\medskip%
   \ifindexmode\write\tabcon{#1 \the\count0}\fi\nobreak}

% TEXT OR MATH MACROS
\let\ep=\epsilon
\let\epsilon=\varepsilon

\def\cite#1{[{\ninerm Ne~Ch.~#1}]}
\def\ii#1.{{\it#1}}

\let\phi=\varphi

\def\N{\hbox{\bbb N}}

\let\hyphen=\-
\def\-{\relax\hyphen}
\def\ro#1{{\rm #1}}
\def\T{\hbox{\eightsans T}}
\def\F{\hbox{\eightsans F}}
\def\FI{\hbox{$\cal F$}}
\def\bul{\hfill\hbox{\tenmsa\char3}}
\def\Th{\medskip\par T{\eightrm HEOREM}}
\def\med#1.#2\par{\smallskip\par\noindent #1.\quad$#2$\par}

\def\pf{\par\smallskip\leavevmode\hskip-10pt\hbox{\nineit Proof. }}

\let\Neg=\neg
\def\neg{\Neg\;}
\def\noin{\smallskip\par\noindent}

\def\nn#1.{$(#1)$}
\def\noj{\medskip\par\quad}
\def\nok// #1.{\par\smallskip\par\noindent(#1)\quad}
\def\Nok/ #1.{\par\smallskip\par\noindent(#1)\quad}

\def\pphi|#1|{\phi(#1)}
\def\pphiz|#1|{\phi^0(#1)}
\def\pphio|#1|{\phi^1(#1)}
\def\pphit|#1|{\phi^2(#1)}
\def\ppphiz|#1|{\phi^{(0)}(#1)}
\def\ppphio|#1|{\phi^{(1)}(#1)}
\def\ppphit|#1|{\phi^{(2)}(#1)}
\long\def\omit#1{}

%NEW

\let\Norm=\|

\def\|{{\scriptstyle\le}}
\def\for{\quad\hbox{for}\quad}
\def\0{\hbox{\rm o}}
\def\1{\hbox{\i}}
\def\2{\hbox{\eighttt,}}
\def\3{\hbox{\eighttt:}}

\def\ari#1{\ul{\hskip.2pt#1\hskip.2pt}\ur}
\def\Open{\vphantom|^{\scriptscriptstyle\{}}
\def\Clos{\vphantom|^{\scriptscriptstyle\}}}
\def\aari#1{\Open{\hskip.2pt#1\hskip.2pt}\Clos}

\let\lbracket=[
\let\rbracket=]
\def\[{\,\big\lbracket\,}
\let\lbrace=\{
\let\rbrace=\}

\def\Q{\hbox{\ninesans Q}}
\def\K{\hbox{$\cal K$}}

%\yy\def\pa{www.math.princeton.edu/~nelson/books/pa.pdf}\zz
\def\PA{\hbox{\ninesans P}}

\def\F{\hbox{\ninesans F}}

\def\sPA{\hbox{\eightsans P}}
\def\ZFC{\hbox{\ninesans ZFC}}
\def\free{\hbox{\sl free }}

\def\co#1{#1^{\hbox{\rm c}}}
\def\Qu[#1]{\langle#1\rangle}

%%%%%%%%%%%%%%%%%%%%%%%%%%%%%%%%%%%%%%%

\def\pr{\'P'}
\def\bw{{^{\scriptscriptstyle<\!\!-}}}
\def\fw{{^{\scriptscriptstyle-\!\!>}}}
\def\tauu_#1{\tau_{\scriptscriptstyle#1}}

\let\lblank=\ %

\def\sind{sind_\phi}

\def\z{\'S'^0}
\def\o{\'S'^1}
\def\SS{\hbox{\ninesans S}}
\def\SSs{\hbox{\eightsans S}}

%% file: BussTao_Afterword_arXiv.tex
\ifx\wrapperUsed\nil
\font\sambigrm=cmr10 scaled \magstep2
\fi
\font\ninett=cmtt9
\font\ninerm=cmr9
\font\ninebf=cmbx9
\font\nineit=cmti9
\font\ninei=cmmi9
\font\ninett=cmtt9
\font\eighti=cmmi8
\font\eightrm=cmr8
\font\eightbf=cmbx8
\font\eightit=cmti8
\font\eightsl=cmsl8
\font\eighttt=cmtt8
\font\eightsy=cmsy8

\font\ninesy=cmsy9
\font\eightsy=cmsy8
\font\nineit=cmti9
\font\eightit=cmti8
\font\sixrm=cmr6
\font\sixbf=cmbx6
\font\sixi=cmmi6
\font\sixsy=cmsy6
\def\ninepoint{\def\rm{\fam0\ninerm}%
   \textfont0=\ninerm \scriptfont0=\sixrm \scriptscriptfont0=\fiverm
   \textfont1=\ninei \scriptfont1=\sixi \scriptscriptfont1=\fivei
   \textfont2=\ninesy \scriptfont2=\sixsy \scriptscriptfont2=\fivesy
   \textfont3=\tenex \scriptfont3=\tenex \scriptscriptfont3=\tenex
   \textfont\itfam=\nineit \def\it{\fam\itfam\nineit}%
   \textfont\bffam=\ninebf \scriptfont\bffam=\sixbf
    \scriptscriptfont\bffam=\fivebf \def\bf{\fam\bffam\ninebf}%
   \normalbaselineskip=11pt
   \setbox\strutbox=\hbox{\vrule height8pt depth 3pt width0pt}%
   \let\big=\ninebig \normalbaselines\rm}

\ifx\wrapperUsed\relax
\magnification=\magstep1
\fi
\hsize=5.0in
\hoffset=0.25in

\def\bibindent{\vskip0.5em\par\noindent\hangindent=1em\hangafter=1}
\def\biblabel#1{#1\hskip4pt\ignorespaces}
\def\BussNelsonPred{[1]}
\def\FarisNelson{[2]}
\def\KritchmanRazsurprise{[3]}
\def\NelsonIST{[4]}
\def\NelsonPredArith{[5]}
\def\YesseninVolpin{[6]}

\def\IDeltazero{{\rm I} \Delta_0}

\
\vskip1cm

\centerline{
\sambigrm  Afterword on two works by Ed Nelson}

\vskip1em

\centerline{
Sam Buss
~~ and ~~ Terence Tao
}

\vskip 2em
Two of Ed Nelson's unfinished papers, {\it Elements} (dated March 12, 2013)
and {\it Inconsistency of Primitive Recursive Arithmetic} (undated, also known
as the {\it Balrog} paper)
have the goal of proving the inconsistency of
number theory.
As Nelson writes in {\it Elements}, ``The aim of this work is to show that
contemporary mathematics, including Peano arithmetic, is inconsistent ...''.

Neither of these papers have been circulated before in their
current forms.  An earlier version of {\it Elements} was circulated
in 2011, but was found to have problems in its treatment of
proofs generated by Chaitin machines.  The new 2013 version
uses a similar approach, but gives a much more detailed
explanation of the planned proof, and it
handles Chaitin machines differently so as to address
the earlier problems.

Nelson's remarkable program to establish the inconsistency
of Peano arithmetic was intertwined with his development
of Internal Set Theory~\NelsonIST{}
and especially Predicative Arithmetic~\NelsonPredArith{}.
Predicative Arithmetic is a constructive fragment
of arithmetic, and Nelson's development of Predicative
Arithmetic was inspired in part
by Yessenin-Volpin's ultra-intuitionistic
set theory~\YesseninVolpin{}.
The mathematical content of Predicate Arithmetic
is closely tied to theories of bounded
arithmetic such as $\IDeltazero$, $\IDeltazero+\Omega_1$,
${\rm S} ^i_2$ and~${\rm T}^i_2$.
Indeed, Nelson~\NelsonPredArith{}
independently discovered some of the important
tools for bounded arithmetic, including the technique of
speeding up induction on cuts and the local interpretability
of predicative arithmetic and bounded arithmetic in
Robinson's theory~$Q$.  Nelson's Predicative Arithmetic
was also influential for the definition by one of us (Buss)
of the theories ${\rm S}^i_2$ and~${\rm T}^i_2$ of
bounded arithmetic, including notably the use of
Nelson's smash function.

The {\it Elements} manuscript gives a detailed,
high-level outline of Nelson's plan for a proof
of the inconsistency of Peano arithmetic (and
primitive recursive arithmetic).  One
of the principal tools is a novel use
of a recent proof
by Kritchman and Raz~\KritchmanRazsurprise{}
of G\"odel's second incompleteness theorem
based on the ``surprise examination''.
Nelson also uses Kolmogorov complexity and
techniques from cut-elimination.  The detailed plan
of the inconsistency proof
is outlined as Steps 1-17 in section~9 near the
end of {\it Elements}.  The plan first discusses a
system~$S^*$ which is a predicative theory including
bounded induction and which is strong enough to express
concepts about metamathetical concepts,
Chaitin machine computation, and Kolmogorov complexity.
Step~7 introduces a finitary theory~$\cal F$; the details
of the system~$\cal F$ are not fully specified, but it needs
to be able to formalize cut-elimination or normalization.
Thus it seems that $\cal F$ can be taken to be, for instance,
$\IDeltazero+{\rm superexp}$ or ${\rm I} \Sigma_1$.
The heart of the argument is reached in Step~16.
Unfortunately, the argument becomes very uncertain here.
Nelson argues that $S^*$ disproves a sequence of statements:
first $A_{\kappa, 0}$, then $A_{\kappa,1}$, etc.,
up through $A_{\kappa,I}$.
The base case that $S^*$ disproves $A_{\kappa,0}$ is fine,
but the later stages are unclear.
It seems that the disproof of $A_{\kappa,\delta}$
requires an assumption that $S^*$ is consistent.
The reason for this is that the Chaitin machine cannot be
given the value of~$\delta$ as an input since the
Kolmogorov complexity of~$\delta$ may not be sufficiently
below that of~$\kappa$. The only alternative to
explicitly specifying~$\delta$ that we can think of,
is for the Chaitin machine to first search
for the $S^*$~proof that ``$\lnot A_{\kappa,\delta+1}$''
and then also wait until $\delta$ many strings are
found to have Chaitin complexity less than~$\kappa$.
This however assumes that $S^*$ is consistent.
Of course, $S^*$ does not prove its own consistency.
Perhaps Nelson had a different argument in mind,
but this is our best attempt to flesh out his arguments.
At any rate,
Nelson was apparently aware of the potential problem here,
since he earlier discusses the need for a system to prove
the ``consistency of its own arithmetization''.

In the spirit of a quote by Carl Sagan, ``Extraordinary claims
require extraordinary evidence'', Nelson planned to fulfill his
inconsistency proof by exhibiting a fully formal,
computer-verified derivation of a contradiction.  That is,
he planned not to prove that there is a proof of contradiction,
but to actually exhibit an explicit proof of a contradiction.  The
first steps of this are carried out at the end of {\it Elements},
and it is even further pursued in {\it Balrog}.
The {\it Balrog} manuscript is still incomplete, as only
six sections are complete, and at least ten sections were planned.
The {\it Balrog} manuscript is in essence
a formalization of the ``bootstrapping''
of predicative arithmetic in the spirit
of~\NelsonPredArith{}.  A remarkable
feature of {\it Balrog} is that
proofs of theorems are indicated in a terse fashion that
permits a Perl program, called {\it qea}, to automatically verify the
proofs.  For instance, Theorems {\it 13b.}\
and {\it 13i.} of~{\it Balrog},
and their proofs, are typeset with
the TeX code
\vskip1em

{\ninett
\noindent
\string\" \string\t\string/\string/13b. 0 + 0 = 0 + 0 \string\"
\vskip1em

\noindent
\string\sam13b.

\noindent
\string\"\string\p\string/13b.

\noindent
\string/\string\'H' \string\ {}

\noindent
\string/5 ; 0 + 0 \string\ {}

\noindent
\string\"
\vskip1em

\noindent
\string\" \string\t\string/13i. x + 0 = 0 + x \string\imp{} \string\'S' x + 0 = 0 + \string\'S' x \string\"
\vskip1em

\noindent
\string\sam13i.

\noindent
\string\"\string\p\string/13i.

\noindent
\string/\string\'H' : x \string\ {}

\noindent
\string/12 ; \string\'S' x \string\ {}

\noindent
\string/12 ; 0 ; x \string\ {}

\noindent
\string/12 ; x \string\ {}

\noindent
\string\"

}
\vskip 1em

These indicate that {\it 13b.}\ is proved by
substituting 0+0 for~$x$ in axiom {\it a5.}, and that
{\it 13i.}\ is proved by using
definition {\it r12.}\ three times,
first substituting
$Sx$ for~$x$, then $0$ and~$x$
for $x$ and~$y$, and finally
$x$ for~$x$.  After these substitutions, the
desired conclusions follow propositionally 
from equality axioms.  The
{\it qea} system then automatically generated
an expanded proof; the expanded proof was produced
as a TeX file, and automatically converted to PDF.
An example is shown in an appendix to the
{\it Balrog} manuscript posted to the {\it arXiv}.

We of course believe that Peano
arithmetic is consistent; thus we do not expect
that Nelson's project can be completed according to
his plans.  Nonetheless, there is much new
in his papers that is of potential mathematical,
philosophical and computational interest.  For this reason,
they are being posted to the {\it arXiv}.  Two aspects
of these papers seem particularly useful.  The first aspect is
the novel use of the ``surprise examination'' and Kolmogorov
complexity; there is some possibility that similar techniques might
lead to new separation results for fragments of
arithmetic.  The second aspect is Nelson's automatic proof-checking
via TeX and {\it qea}.  This is highly
interesting and provides
a novel method of integrating human-readable
proofs with computer verification of proofs.
\medskip

The reader interested in further discussion of
Nelson's Predicative Arithmetic can
consult the mostly-survey article~\BussNelsonPred{}.
The volume~\FarisNelson{} contains papers about
many other aspects of Nelson's wide-ranging research.
Other works by Nelson are available at
{\tt math.princeton.edu\slash$\sim$nelson}, including a number
of philosophical works.

\vskip0.5em

{\ninepoint
\bibindent
\biblabel{\BussNelsonPred}
{S.~R.\ Buss}, {\it Nelson's work on logic and foundations and other
  reflections on foundations of mathematics}, in Diffusion, Quantum Theory, and
  Radically Elementary Mathematics, Princeton University Press, 2006,
  pp.~183--208. Edited by W.\ Faris.
\bibindent
\biblabel{\FarisNelson}
{W.~G.\ Faris}, ed., {\it Diffusion, Quantum Theory, and Radically
  Elementary Mathematics}, Mathematical Notes, \#47, Princeton University
  Press, 2006.
\bibindent
\biblabel{\KritchmanRazsurprise}
{S.~Kritchman and R.~Raz}, {\it The surprise examination and the second
  incompleteness theorem}, Notices of the American Mathematical Society, 57
  (2010), pp.~1454--1458.
\bibindent
\biblabel{\NelsonIST}
{E.~Nelson}, {\it Internal set theory: A new approach to nonstandard
  analysis}, Bulletin of the American Mathematical Society, 83 (1977),
  pp.~1165--1198.
\bibindent
\biblabel{\NelsonPredArith}
{E.~Nelson}, {\it {P}redicative {A}rithmetic}, Princeton University Press,
  1986.
\bibindent
\biblabel{\YesseninVolpin}
{A.~S.\ Yessenin-Volpin}, {\it The ultra-intuitionistic criticism and the
  antitraditional program for foundations of mathematics}, in Intuitionism and
  Proof Theory, A.~Kino, J.~Myhill, and R.~E. Vesley, eds., North-Holland,
  1970, pp.~1--45.
}

\bigskip
\hfill
{Sam Buss}
\vskip 0.5ex

\hfill
{Terence Tao}
\vskip 2ex

\hfill {\it
August 31, 2015}

\ifx\wrapperUsed\nil
\def\thisEnd{\end}
\else
\def\thisEnd{\relax}
\fi

\thisEnd